\documentclass[11pt,a4paper,oneside]{article}%
\usepackage{amsfonts}
\usepackage{graphicx}
\usepackage{amsmath,amssymb}
\usepackage{color}
\usepackage{amssymb}
\usepackage{amsmath}
\usepackage[T1]{fontenc}
\usepackage{amsthm}
\usepackage[margin=1in]{geometry}
\usepackage{amssymb}
\usepackage{amsfonts}
\usepackage{graphicx}
\usepackage{amsmath}%
\providecommand{\U}[1]{\protect\rule{.1in}{.1in}}
\newtheorem{theorem}{Theorem}

\newtheorem{condition}[theorem]{Condition}

\newtheorem{corollary}[theorem]{Corollary}

\newtheorem{lemma}[theorem]{Lemma}

\newtheorem{proposition}[theorem]{Proposition}
\newtheorem{remark}[theorem]{Remark}

\begin{document}

\begin{center}
{\Large \textbf{On the empirical spectral distribution for matrices with long
memory and independent rows }}

\bigskip

F. Merlev\`{e}de and M. Peligrad\footnote{Supported in part by a Charles
Phelps Taft Memorial Fund grant, and the NSF grant DMS-1208237.}

\bigskip

\end{center}

Universit\'{e} Paris Est, LAMA (UMR 8050), UPEM, CNRS, UPEC. 

Email: florence.merlevede@u-pem.fr \bigskip

Department of Mathematical Sciences, University of Cincinnati, PO Box 210025,
Cincinnati, Oh 45221-0025, USA.

Email: peligrm@ucmail.uc.edu


\bigskip

\begin{center}
\bigskip\textbf{Abstract}
\end{center}

In this paper we show that the empirical eigenvalue distribution of any sample
covariance matrix generated by independent copies of a stationary regular
sequence has a limiting distribution depending only on the spectral density of
the sequence. We characterize this limit in terms of Stieltjes transform via a
certain simple equation. No rate of convergence to zero of the covariances is
imposed. If the entries of the stationary sequence are functions of
independent random variables the result holds without any other additional assumptions.

As a method of proof, we study the empirical eigenvalue distribution for a
symmetric matrix with independent rows below the diagonal; the entries satisfy
a Lindeberg-type condition along with mixingale-type conditions without rates.
In this nonstationary setting we point out a property of universality, meaning
that, for large matrix size, the empirical eigenvalue distribution depends
only on the covariance structure of the sequence and is independent on the
distribution leading to it. These results have interest in themselves,
allowing to study symmetric random matrices generated by random processes with
both short and long memory.

\section{Introduction and the main Result.}

\ Due to the fact that random matrices appear in many applied fields, their
empirical spectral distribution is a subject of intense research. Earlier
works, pioneered by the celebrated paper by Wigner (1958), deal with symmetric
matrices having independent entries below the diagonal. Only in the last two
decades there has been an effort to weaken the hypotheses of independence and
various forms of weak dependence have been considered. The progress was in
general achieved first for Gaussian random matrices. For this case the joint
distribution of eigenvalues is tractable. Among the papers for symmetric
Gaussian matrices with correlated entries we mention the works of Khorunzhy
and \ Pastur (1994), Boutet de Monvel \textit{et al.} (1996), Boutet de Monvel
and Khorunzhy (1999), Chakrabarty \textit{et al.} (2014).

Our paper is essentially motivated by the study of large sample covariance
matrices, which is a very important topic in multivariate analysis. The
spectral analysis of large-dimensional sample covariance matrices has been
actively studied starting with the work of \ Mar\u{c}enko and Pastur (1967).
Extensions can be found in the works of Wachter (1978), Yin (1986),
Silverstein (1995), Silverstein and Bai (1995), Hachem \textit{et al. }(2005),
Bai and Zhou (2008), Adamczak (2011), Pfaffel and Schlemm (2011), Yao (2012),
Banna and Merlev\`{e}de (2013).

In this paper, in Theorem \ref{corgram}, we find the limiting empirical
eigenvalue distribution for the sample covariance matrix of a stationary
process which is regular. The regularity is an ergodic-type property. Our
result shows that the limit can be obtained much beyond the situation of
weakly dependent case which corresponds to continuous and bounded spectral
densities, short range dependence and absolutely summable covariances. It
applies to long range dependent stationary stochastic processes and sheds
light on the theory of sample covariance random matrices, which is important
in large sample statistics for stochastic processes. We show that the limit of
the empirical spectral distribution exists and we also characterize the limit
in terms of its Stieltjes transform, even for the case when the spectral
density is not continuous or even square integrable. The previous works on
short memory processes heavily relied upon the limits of Toeplitz matrices
induced by the \ covariance structure. In our case, the covariance matrices
fall outside the range of the celebrated Szeg\"{o}-Trotter theorem for
Toeplitz matrices, which is restricted to square summable entries, therefore
to square integrable spectral densities. We showed that the spectral density
of the underlying stationary process is a key factor to describe the limit in
all the situations.

Furthermore, the technical theorems leading to Theorem \ref{corgram} are also
important. They reduce the study of the empirical spectral distribution of
symmetric matrices with independent regular rows, below diagonal, to the study
of the sequence of the expected value of Stieltjes transforms associated to a
Gaussian matrix with the same covariance structure. These results are set in
the non-stationary case for variables satisfying a certain Lindeberg
condition. Their proofs are complicated by the fact that our intention was to
avoid the use of rates of decay of the covariances.

In order to stress the importance of our results we include several
applications to regular processes, functions of i.i.d., and linear processes
with martingale differences innovations. As we shall see, Theorem
\ref{corgram} applies to large sample covariance matrices constructed from
independent copies of any stationary process whose entries are functions of
i.i.d. which are centered and has finite second moments. In particular the
theorem applies to any causal linear process with square summable coefficients
and i.i.d. innovations as soon as the process exists in ${\mathbb{L}}^{2}$, so
it could have long memory.

Our proofs are a blend of probabilistic techniques for dependent structures
such as the big and small block argument and martingale approximations,
properties of Gaussian processes, and algebraic and Fourier analysis tools.
Because our variables are correlated the method of proof is based on the
Stieltjes transform, which is well adapted to handle dependent entries. The
Stieltjes transform is also useful to characterize the limit.

Here are some notations used all along the paper. The notation $[x]$ is used
to denote the integer part of a real $x$. The notation $\mathbf{0}_{p}$ means
a row vector of size $p$ with components equal to zero. When no confusion is
possible concerning the size of a null vector ${\mathbf{0}}$ we will omit the
index of its size. For a matrix $A$, we denote by $A^{T}$ its transpose
matrix, by $\mathrm{Tr}(A)$ its trace. We shall also use the notation $\Vert
X\Vert_{r}$ for the $\mathbb{L}^{r}$-norm ($r\geq1$) of a real valued random
variable $X$.

For any sequence of square matrices $A_{n}$ of order $n$ with only real
eigenvalues $\lambda_{1,n}\leq\dots\leq\lambda_{n,n}$, the spectral
distribution function is defined by
\[
\ F^{A_{n}}(x)=\frac{1}{n}\sum_{k=1}^{n}I(\lambda_{k,n}\leq x)\,,
\]
where $I(B)$ denotes the indicator of an event $B$. The general problem is to
find a distribution function $F$ such that $F^{A_{n}}\rightarrow F$ at all
points of continuity of $F,$ or equivalently $d(F^{A_{n}},F)\rightarrow0,$
where the L\'{e}vy distance between two distribution functions $F$ and $G$ is
defined by
\[
d(F,G)=\inf\{\varepsilon>0\ :\ F(x-\varepsilon)-\varepsilon\leq G(x)\leq
F(x+\varepsilon)+\varepsilon\}\,.
\]
The Stieltjes transform of $F^{A_{n}}$ is given by
\[
S^{A_{n}}(z)=\int\frac{1}{x-z}dF^{A_{n}}(x)=\frac{1}{n}\mathrm{Tr}%
(A_{n}-z{\mathbf{I}}_{n})^{-1}\,,
\]
where $z=u+\mathrm{i}v\in\mathbb{C}^{+}$ (the set of complex numbers with
positive imaginary part), and $\mathbf{I}_{n}$ is the identity matrix of order
$n$.$\ $\ It is well-know that $\lim_{n\rightarrow\infty}d(F^{A_{n}},F)=0$ if
and only if for all $z\in\mathbb{C}^{+}$, $S_{A_{n}}(z)\rightarrow S_{F}(z).$
We can also see, for instance, in Proposition 2.1 in Bobkov \textit{et al.}
(2010), that the estimate of the L\'{e}vy distance between empirical spectral
distribution functions associated with two matrices can be also given in terms
of their Stieltjes transforms.

Let $N$ and $p$ be two positive integers and consider the $N\times p$ matrix
\begin{equation}
{\mathcal{X}}_{N,p}=\big (X_{ij}\big )_{1\leq i\leq N,1\leq j\leq p}\,,
\label{defcalX}%
\end{equation}
where $X_{ij}$'s are real-valued random variables. Define now the symmetric
matrix ${\mathbb{B}}_{N}$ of order $p$ by
\begin{equation}
{\mathbb{B}}_{N}=\frac{1}{N}{\mathcal{X}}_{N,p}^{T}{\mathcal{X}}_{N,p}\,.
\label{defBb}%
\end{equation}
The matrix ${\mathbb{B}}_{N}$ is usually referred to as the sample covariance
matrix associated with the process $(X_{\mathbf{u}})_{\mathbf{u}\in
{\mathbb{Z}}^{2}}$. It is also known under the name of Gram random matrix.

In Theorem \ref{corgram} below, we consider $N$ independent copies
$(X_{ij})_{j\in{\mathbb{Z}}}$, $i=1,\dots,N$ of a stationary sequence
$(X_{i})_{i\in{\mathbb{Z}}}$ of real-valued random variables in ${\mathbb{L}%
}^{2}$ and give sufficient conditions to characterize the limiting
distribution of $F^{{\mathbb{B}}_{N}}$ when $p/N\rightarrow c\in(0,\infty)$.
Relevant to this characterization is the notion of spectral distribution
function induced by the covariances of $(X_{i})_{i\in{\mathbb{Z}}}$. By
Herglotz's Theorem (see e.g. Brockwell and Davis \cite{BD}), there exists a
non-decreasing function $G$ (the spectral distribution function) on $[-\pi
,\pi]$ such that, for all $j\in\mathbb{Z}$, $\mathrm{Cov}(X_{0},X_{j}%
)=\int_{-\pi}^{\pi}\exp({\mathrm{i}}j\theta)dG(\theta)$. If $G$ is absolutely
continuous with respect to the normalized Lebesgue measure $\lambda$ on
$[-\pi,\pi]$, then the Radon-Nikodym derivative $f$ of $G$ with respect to the
Lebesgue measure is called the spectral density, it is a nonnegative, even and
integrable function on $[-\pi,\pi]$ which satisfies
\[
\mathrm{Cov}(X_{0},X_{j})=\int_{-\pi}^{\pi}\exp({\mathrm{i}}j\theta
)f(\theta)d\theta\,,\ \ j\in\mathbb{Z}\,.
\]
We shall introduce the following regularity conditions. Define the left tail
sigma field of $(X_{i})_{i\in{\mathbb{Z}}}$ by ${\mathcal{G}}_{-\infty
}=\bigcap_{k\in{\mathbb{Z}}}{\mathcal{G}}_{k}$ where ${\mathcal{G}}_{k}%
=\sigma(X_{j},j\leq k)$
\begin{equation}
\mathbb{E}(X_{0}|\mathcal{G}_{-\infty})=0\text{ a.s.} \label{regular 1}%
\end{equation}
and for every integer $k$
\begin{equation}
\mathbb{E}(X_{0}X_{k}|\mathcal{G}_{-\infty})=\mathbb{E}(X_{0}X_{k})\text{
a.s.} \label{regular 2}%
\end{equation}
We point out that if \eqref{regular 1} holds, then the process $(X_{k}%
)_{k\in{\mathbb{Z}}}$ is purely non deterministic. Hence, by a result of
Szeg\"{o} (see for instance \cite[Theorem 3]{Bingh}) if \eqref{regular 1}
holds, the spectral density $f$ of $(X_{k})_{k\in{\mathbb{Z}}}$ exists and if
$X_{0}$ is non degenerate,
\[
\int_{-\pi}^{\pi}\log f(t)~dt>-\infty\,;
\]
in particular, $f$ cannot vanish on a set of positive measure.

\begin{theorem}
\label{corgram} Consider $N$ independent copies $(X_{ij})_{j\in{\mathbb{Z}}}$,
$i=1,\dots,N$ of a stationary sequence $(X_{i})_{i\in{\mathbb{Z}}}$ of
real-valued random variables centered and in ${\mathbb{L}}^{2}$ and that
satisfies the conditions \eqref{regular 1} and \eqref{regular 2}. Assume
$p/N\rightarrow c\in(0,\infty)$. Then there is a nonrandom probability
distribution $F$ such that $d(F^{{\mathbb{B}}_{N}},F)\rightarrow0$ a.s.
Furthermore, the Stieltjes transform $S=S(z),$ $z\in\mathbb{C}^{+}$, of $F$ is
determined by the equation
\begin{equation}
z=-\frac{1}{{\underline{S}}}+\frac{c}{2\pi}\int_{-\pi}^{\pi}\frac
{1}{{\underline{S}}+(2\pi f(\lambda))^{-1}}d\lambda\,, \label{equationlimit}%
\end{equation}
where ${\underline{S}}:=-(1-c)/z+cS$ and $f(\cdot)$ is the spectral density of
$(X_{k})_{k\in\mathbb{Z}}$.
\end{theorem}

\begin{remark}
As a matter of fact, we can relax the stationarity to stationarity in
${\mathbb{L}}^{2}$. More precisely, the conclusion of Theorem \ref{corgram}
applies for Gram matrices constructed from a process $(X_{\mathbf{u}%
})_{\mathbf{u}\in{\mathbb{Z}}^{2}}$ satisfying the conditions of Theorem
\ref{mainindependentrows} below if we assume in addition that for any
$i,k,\ell$ in ${\mathbb{Z}}$
\[
\mathrm{Cov}(X_{ik},X_{i\ell})=\mathrm{Cov}(X_{0k},X_{0\ell})=\mathrm{Cov}%
(X_{00},X_{0,\ell-k})\,.
\]
In this case, $f(\cdot)$ is the spectral density of $(X_{0k})_{k\in\mathbb{Z}%
}$.
\end{remark}

Note that if ${\mathcal{G}}_{-\infty}$ is trivial then the conditions
\eqref{regular 1} and \eqref{regular 2} hold. Therefore we can immediately
formulate the following corollary to Theorem \ref{corgram}:

\begin{corollary}
Consider $N$ independent copies $(X_{ij})_{j\in{\mathbb{Z}}}$, $i=1,\dots,N$
of a stationary sequence $(X_{i})_{i\in{\mathbb{Z}}}$ of real-valued random
variables centered and in ${\mathbb{L}}^{2}$ with trivial left tail sigma
field, ${\mathcal{G}}_{-\infty}$. Then the conclusion of Theorem \ref{corgram} holds.
\end{corollary}

\section{Some technical results for symmetric matrices}

A key step in the proof of Theorem \ref{corgram} is to show that the study of
the limiting spectral distribution function of ${\mathbb{B}}_{N}$ can be
reduced to studying the same problem as for a Gaussian matrix with the same
covariance structure. This step will be achieved with the help of some
preliminary technical results concerning symmetric matrices with independent
rows below the diagonal. These technical results have interest in themselves
since they show that, for symmetric matrices with independent rows below the
diagonal, very simple regularity conditions on the entries of each row allow
to reduce the study of their limiting spectral distribution function to the
one of a symmetric Gaussian matrix with the same covariance structure. In
particular, this applies when the rows, below the diagonal, are independent
and generated by the same stationary sequence provided it is regular, i.e. has
a trivial left tail sigma-field.

\medskip

To state the results of this section, let us introduce some notations. Let
$(X_{\mathbf{u}})_{\mathbf{u}\in\mathbb{N}^{2}}$ be real-valued random
variables on a probability space $(\Omega,\mathcal{F},\mathbb{P})$. In what
follows, we consider the symmetric $n\times n$ random matrix $\mathbf{X}_{n}$
defined as follows: for any $i$ and $j$ in $\{1,\dots,n\},$
\begin{align}
(\mathbf{X}_{n})_{ij}  &  =X_{ij}\,\text{ for }i\geq j\ \text{ and
}\label{defX}\\
(\mathbf{X}_{n})_{ij}  &  =X_{ji}\,\text{ for }i<j\,.\nonumber
\end{align}
Define
\begin{equation}
{\mathbb{X}}_{n}:=\frac{1}{n^{1/2}}\mathbf{X}_{n}\,, \label{defW}%
\end{equation}
and set
\[
L(A)=\frac{1}{n^{2}}\sum_{i=1}^{n}\sum_{j=1}^{i}{\mathbb{E}}(X_{ij}%
^{2}I(|X_{ij}|>A))\,,
\]
where $A$ is a positive number.

We shall introduce now a Lindeberg's type condition:

\begin{condition}
\label{condition 1}(1) ${\mathbb{E}}(X_{\mathbf{u}})=0$ for all $\mathbf{u}%
\in\mathbb{N}^{2}.$ \newline(2) There is $\sigma>0$ such that $\sup
_{\mathbf{u}\in\mathbb{N}^{2}}\Vert X_{\mathbf{u}}\Vert_{2}\leq\sigma.$
\newline(3) \ For every $\varepsilon>0$ we have $L(\varepsilon n^{1/2}%
)\rightarrow0$ as $\rightarrow\infty.$
\end{condition}

Clearly the items (2) and (3) of this condition are satisfied as soon as the
family $(X_{\mathbf{u}}^{2})$ is uniformly integrable or the random field is stationary.

Next result, in the nonstationary setting, shows that two mild regularity-like
conditions without rates, are sufficient for reducing the study of the
limiting spectral distribution of a symmetric matrix with independent rows
below the diagonal to the corresponding problem for a Gaussian matrix having
the same covariance structure. This result indicates that for large matrix
size, the empirical distribution of the eigenvalues is universal, in the sense
that it is determined only by the covariance structure of the process.

\begin{theorem}
\label{mainindependentrows} Assume that Condition \ref{condition 1} is
satisfied and in addition that the random vectors $(R_{i})_{i\geq1},$ where
$R_{i}=(X_{ij})_{j\in\mathbb{N}},$ are mutually independent. For any $i\geq1$
fixed, let ${\mathcal{G}}_{ik}=\sigma(X_{ij},1\leq j\leq k)$ and, by
convention, for $k\leq0,$ ${\mathcal{G}}_{ik}=\{\emptyset,\Omega\}$. Then,
under the following two additional assumptions:
\begin{equation}
\eta_{m}=\sup_{i \geq j\geq m}\Vert\mathbb{E}(X_{ij}|\mathcal{G}_{i,j-m}%
)\Vert_{2}\rightarrow0 \label{mart-K-indrow}%
\end{equation}
and
\begin{equation}
\gamma_{m}=\sup_{i \geq\ell\geq k\geq m}\Vert\mathbb{E}(X_{ik}X_{i\ell
}|\mathcal{G}_{i,k-m})-\mathbb{E}(X_{ik}X_{i\ell})\Vert_{1}\rightarrow0\,,
\label{mixing-indrow}%
\end{equation}
the following convergence holds: for all $z\in\mathbb{C}^{+}$,
\begin{equation}
S^{\mathbb{X}_{n}}(z)-\mathbb{E}S^{\mathbb{Y}_{n}}(z)\rightarrow0\text{ almost
surely, as }n\rightarrow\infty, \label{conclusionas}%
\end{equation}
where ${\mathbb{X}}_{n}$ is defined by (\ref{defW}) and $\mathbb{Y}%
_{n}=\mathbf{Y}_{n}/\sqrt{n}$, $\mathbf{Y}_{n}$ being the symmetric matrix
defined as in \eqref{defX} and constructed from a centered real-valued
Gaussian random field $(Y_{\mathbf{u}})_{\mathbf{u}\in{\mathbb{N}}^{2}}$
having the same covariance structure as $(X_{\mathbf{u}})_{\mathbf{u}%
\in{\mathbb{N}}^{2}}.$
\end{theorem}

\begin{remark}
Since $\mathbf{Y}_{n}$ is constructed from a centered real-valued Gaussian
random field $(Y_{\mathbf{u}})_{\mathbf{u}\in{\mathbb{N}}^{2}}$ having the
same covariance structure as $(X_{\mathbf{u}})_{\mathbf{u}\in{\mathbb{N}}^{2}%
}$, we have in particular that the random vectors $(G_{i})_{i\geq1},$ where
$G_{i}=(Y_{ij})_{j\in\mathbb{N}}$, are mutually independent. Therefore
relation (\ref{fact1TH2}) in the proof of Theorem \ref{mainindependentrows}
also holds for $\mathbb{Y}_{n}$. Hence, in addition to the conclusion of
Theorem \ref{mainindependentrows}, we also have
\[
S^{\mathbb{X}_{n}}(z)-S^{\mathbb{Y}_{n}}(z)\rightarrow0\text{ almost surely,
as }n\rightarrow\infty,
\]
provided that $(X_{\mathbf{u}})_{\mathbf{u}\in\mathbb{N}^{2}}$ and
$(Y_{\mathbf{u}})_{\mathbf{u}\in\mathbb{N}^{2}}$ are defined on the same
probability space.
\end{remark}

\begin{remark}
\label{remtriangular} Theorem \ref{mainindependentrows} also holds if we allow
the random variables $X_{ij}$ to depend on the matrix size $n$. In this
context we write $X_{ij}^{(n)}$ instead of $X_{ij}$, we adapt in a natural way
Condition 1 and we modify conditions \eqref{mart-K-indrow} and
\eqref{mixing-indrow} as follows:
\[
\sup_{n\geq1}\sup_{i\geq j\geq m}\Vert\mathbb{E}(X_{ij}^{(n)}|\mathcal{G}%
_{i,j-m}^{(n)})\Vert_{2}\rightarrow_{m\rightarrow\infty}0
\]
and
\[
\sup_{n\geq1}\sup_{i\geq\ell\geq k\geq m}\Vert\mathbb{E}(X_{ik}^{(n)}X_{i\ell
}^{(n)}|\mathcal{G}_{i,k-m}^{(n)})-\mathbb{E}(X_{ik}^{(n)}X_{i\ell}%
^{(n)})\Vert_{1}\rightarrow_{m\rightarrow\infty}0\,.
\]

\end{remark}

Next corollary applies to stationary sequences and shows that the conclusion
of Theorem \ref{mainindependentrows} holds under simple regularity conditions.

\begin{corollary}
\label{stationary}Let $(X_{ij})_{j\in{\mathbb{Z}}}$, $i=1,\dots,n$ be $n$
independent copies of a stationary sequence $(X_{k})_{k\in{\mathbb{Z}}}$ of
real-valued random variables which are centered and in ${\mathbb{L}}^{2}$.
Then the conclusion of Theorem \ref{mainindependentrows} holds under the
regularity conditions \eqref{regular 1} and \eqref{regular 2}.
\end{corollary}

Theorem \ref{mainindependentrows} and its Remark \ref{remtriangular} allow us
to formulate the following result for Gram matrices. It will be a key step in
the proof of Theorem \ref{corgram}.

\begin{theorem}
\label{Gram} Under the conditions of Theorem \ref{mainindependentrows} and if
$p/N\rightarrow c\in(0,\infty)$, the following convergence holds: for all
$z\in\mathbb{C}^{+}$,
\[
S^{\mathbb{B}_{N}}(z)-{\mathbb{E}}S^{\mathbb{H}_{N}}(z)\rightarrow0\text{
almost surely, as }N\rightarrow\infty\, ,
\]
where $\mathbb{B}_{N}$ is defined by (\ref{defBb}) and $\mathbb{H}_{N}$ is a
Gram random matrix associated with a centered real-valued Gaussian process
$(Y_{\mathbf{u}})_{\mathbf{u}\in{\mathbb{Z}}^{2}}$ having the same covariance
structure as $(X_{\mathbf{u}})_{\mathbf{u}\in{\mathbb{Z}}^{2}}$.
\end{theorem}

\section{Examples}

Below we give a few examples of regular processes.

\smallskip\ 

1. \textbf{Functions of i.i.d. random variables.} Let $(\varepsilon
_{\mathbf{u}})_{\mathbf{u}\in{\mathbb{Z}}}$ be i.i.d. and $g:{\mathbb{R}%
}^{{\mathbb{Z}}}\rightarrow{\mathbb{R}}$ be a measurable function such that,
for any $i,j$ in ${\mathbb{Z}}$, $X_{ij}=g(\varepsilon_{ik},k\leq j)$ is well
defined in ${\mathbb{L}}^{2}$ and ${\mathbb{E}}(X_{ij})=0$. These are regular
random fields since each row has a trivial left sigma field. Therefore for
these processes, conditions \eqref{regular 1} and \eqref{regular 2} are
satisfied. Examples include linear processes, functions of linear processes
and iterated random functions (see for instance Wu and Woodroofe (2000), among others).

For example let $X_{ij}=\sum_{k=0}^{\infty}a_{k}\varepsilon_{i,k-j}$, where
$\varepsilon_{ij}$ are i.i.d. with mean $0$ and finite variance, and $a_{k}$
are real coefficients with $\sum_{k=1}^{\infty}a_{k}^{2}<\infty$. In this case
$X_{ij}$ is well-defined, the process is regular, and therefore the conclusion
of Theorem \ref{corgram} holds. The limiting empirical eigenvalue distribution
of Gram matrices associated with linear processes was investigated in several
papers (see for instance \cite{PS}, \cite{Yao} and \cite{ban}) but, all the
previous known results treat only the short memory case meaning that the
$a_{k}$'s are absolutely summable.

\bigskip

As we mentioned before, conditions (\ref{regular 1}) and (\ref{regular 2}) are
satisfied for a stationary sequence if the left tail sigma field
$\mathcal{G}_{-\infty}$ is trivial. Processes with trivial tail sigma field
are called regular (see Chapter 2, Volume 1 in Bradley, 2007). We give next
examples of regular processes.

1. \textbf{Mixing sequences. }The strong mixing coefficient is defined in the
following way:%
\[
\alpha(\mathcal{A},\mathcal{B)=}\sup\{|{\mathbb{P}}(A\cap B)-{\mathbb{P}%
}(A){\mathbb{P}}(B)|:\,A\in\mathcal{A},\,B\in\mathcal{B}\},
\]
where $\mathcal{A}\ $\ and $\mathcal{B}$ are two sigma algebras.

The $\rho-$mixing coefficient, also known as maximal coefficient of
correlation, is defined as
\[
\rho(\mathcal{A},\mathcal{B})=\sup\{\mathrm{Cov}(X,Y)/\Vert X\Vert_{2}\Vert
Y\Vert_{2}:\,X\in\mathbb{L}^{2}(\mathcal{A}),\,Y\in\mathbb{L}^{2}%
(\mathcal{B})\}.
\]
For the stationary sequence of random variables $(X_{k})_{k\in\mathbb{Z}}$,
$\mathcal{F}^{n}$ denotes the $\sigma$--field generated by $X_{i}$ with
indices $i\geq n,$ and $\mathcal{F}_{m}$ denotes the $\sigma$--field generated
by $X_{i}$ with indices $i\leq m$. Then we define the sequences of mixing
coefficients
\[
\alpha_{n}=\alpha(\mathcal{F}_{0},\mathcal{F}^{n}\mathcal{)}\ \text{ and }%
\rho_{n}=\rho(\mathcal{F}_{0},\mathcal{F}^{n}\mathcal{)} \, .
\]

A sequence is called strongly mixing if $\alpha_{n}\rightarrow0$. It is
well-known that for strongly mixing sequences the left tail sigma field is
trivial; see Claim 2.17a in Bradley (2007). Examples of this type include
Harris recurrent Markov chains.

If $\lim_{n\rightarrow\infty}\rho_{n}<1$, then the tail sigma field is also
trivial according to Section 2.5 in Bradley (2005).

Note that our conditions \eqref{mart-K-indrow} and \eqref{mixing-indrow} also
hold without the assumptions of stationarity and of regularity. For instance,
if
\[
\alpha_{2,n}:=\sup_{i\geq1}\sup_{j\geq k}\alpha\big (\sigma(X_{i1}%
,\dots,X_{ik}),\sigma(X_{i,k+n},X_{i,j+n})\big )\rightarrow0\,,
\]
and if the variables are centered and $(X_{\mathbf{u}}^{2})_{{\mathbf{u}}%
\in{\mathbb{Z}}^{2}}$ is uniformly integrable, then \eqref{mart-K-indrow} and
\eqref{mixing-indrow} are satisfied. Note that the condition $\alpha
_{2,n}\rightarrow0$ is not enough for regularity.

\bigskip

For a nonstationary example we shall look at a more general linear process,
based on martingale difference innovations satisfying Lindeberg's condition.

2. \textbf{Linear processes} \textbf{with martingale entries}. Assume that for
any $1\leq j\leq i\leq n$, the $(i,j)^{\text{th}}$ entry of ${\mathbf{X}}_{n}$
is given by a linear process of the form
\begin{equation}
X_{ij}=\sum_{\ell=0}^{\infty}a_{i\ell}d_{i,j-\ell}\,,\label{deflinearprocess}%
\end{equation}
where $(a_{\mathbf{u}})_{{\mathbf{u}}\in{\mathbb{Z}}^{2}}$ is a sequence of
real numbers and $(d_{{\mathbf{u}}})_{{\mathbf{u}}\in{\mathbb{Z}}^{2}}$ is a
sequence of real-valued random variables satisfying the conditions below:

\begin{enumerate}
\item[${\mathbf{A}}_{1}$] $A_{n,i}=\sum_{j=0}^{n}a_{ij}^{2}<\infty$ is
convergent as $n\rightarrow\infty$ \thinspace uniformly in $i\geq1$. \newline

\item[${\mathbf{A}}_{2}$] There is $\sigma>0$ such that $\sup_{{\mathbf{u}}%
\in{\mathbb{Z}}^{2}}\Vert d_{\mathbf{u}}\Vert_{2}<\sigma$ and for every
$\varepsilon>0$,
\[
\frac{1}{n^{2}}\sum_{i=1}^{n}\sum_{j=1}^{i}{\mathbb{E}}(d_{ij}^{2}%
I(|d_{ij}|>\varepsilon\sqrt{n}))\rightarrow0\,\text{ as $n\rightarrow\infty$%
}.
\]

\item[${\mathbf{A}}_{3}$] Setting ${\mathcal{F}}_{ij}=\sigma(d_{ik},k\leq j)$,
${\mathbb{E}}(d_{ij}|{\mathcal{F}}_{i,j-1})=0$ a.s. for any $(i,j)$ in
${\mathbb{Z}}^{2}$ and
\[
\sup_{i\geq1}\sup_{j\geq n}\Vert{\mathbb{E}}(d_{ij}^{2}|{\mathcal{F}}%
_{i,j-n})-{\mathbb{E}}(d_{ij}^{2})\Vert_{1}\rightarrow0\,\text{ as
$n\rightarrow\infty$}.
\]

\end{enumerate}

\begin{corollary}
\label{corlinear1} Assume that $(X_{ij})$ is a linear process as defined in
\eqref{deflinearprocess} such that the conditions ${\mathbf{A}}_{1}$,
${\mathbf{A}}_{2}$ and ${\mathbf{A}}_{3}$ hold. Assume in addition that the
random vectors $(d_{i{\mathbf{.}}})_{i\geq1},$ where $d_{i{\mathbf{.}}%
}=(d_{ij})_{j\in\mathbb{Z}}$, are mutually independent. Then the conclusion of
Theorem \ref{mainindependentrows} hold.
\end{corollary}

The proof of this corollary is based on standards arguments for martingales
and is left to the reader.

\section{Proofs}

\subsection{Preparatory materials}

In this section, we collect several results useful for our proofs.

The first result we mention is Lemma 2.1 in G\"{o}tze \textit{et al.} (2012)
that allows to compare the difference between two Stieltjes transforms.

\begin{lemma}
\label{lmagotze} Let ${\mathbf{A}}$ and ${\mathbf{B}}$ be two symmetric
$n\times n$ matrices with real entries. Then, for any $z=x+{\mathrm{i}}%
y\in{\mathbb{C}}\backslash{\mathbb{R}}$,
\[
|S_{{\mathbf{A}}}(z)-S_{{\mathbf{B}}}(z)|\leq\frac{1}{y^{2}\sqrt{n}%
}|\mathrm{Tr}({\mathbf{A}}-{\mathbf{B}})|^{1/2}\,.
\]

\end{lemma}

Relevant to the proof of Theorem \ref{corgram} is the following lemma which
gives an estimate of the L\'{e}vy distance between two distribution functions
of eigenvalues (see Corollary A.42 in Bai and Silverstein (2010)).

\begin{lemma}
\label{lmagotze2} Let ${\mathbf{A}}$ and ${\mathbf{B}}$ be two $n\times p$
matrices with real entries, and $d$ be the L\'{e}vy distance. Then, for any
$z=x+{\mathrm{i}}y\in{\mathbb{C}}\backslash{\mathbb{R}}$,
\[
d^{2}(F_{{\mathbf{A}}{\mathbf{A}}^{T}},F_{{\mathbf{B}}{\mathbf{B}}^{T}}%
)\leq\frac{\sqrt{2}}{n}[\mathrm{Tr}({\mathbf{A}}{\mathbf{A}}^{T}+{\mathbf{BB}%
}^{T})\mathrm{Tr}(({\mathbf{A}}-{\mathbf{B}})({\mathbf{A}}-{\mathbf{B}}%
)^{T})]^{1/2}\,.
\]

\end{lemma}

All along the proofs, we shall use the fact that the Stieltjes transform of
the spectral measure is a smooth function of the matrix entries. To formalize
things in a way that is suitable for our purpose, we shall adopt the same
notations as in Chatterjee (2006) and introduce the following map $A$ which
"constructs" Wigner-type matrices. Let $N=n(n+1)/2$ and write elements of
${\mathbb{R}}^{N}$ as $\mathbf{x}=(x_{ij})_{1\leq j\leq i\leq n}$. For any
$\mathbf{x}$ in ${\mathbb{R}}^{N}$, let $A(\mathbf{x})$ be the matrix defined
by
\begin{equation}
(A(\mathbf{x}))_{ij}=\left\{
\begin{array}
[c]{ll}%
\frac{1}{\sqrt{n}}x_{ij}\  & i\geq j\\
\frac{1}{\sqrt{n}}x_{ji}\  & i<j\,.
\end{array}
\right.  \label{defofA}%
\end{equation}
Let $z\in\mathbb{C}^{+}$ and $s_{n}:=s_{n,z}$ be the function defined from
$\mathbb{R}^{N}$ to $\mathbb{C}$ by
\begin{equation}
s_{n}(\mathbf{x})=\frac{1}{n}\mathrm{Tr}(A(\mathbf{x})-z{\mathbf{{I}}_{n}%
})^{-1}\,, \label{deffA}%
\end{equation}
where $\mathbf{{I}}_{n}$ is the identity matrix of order $n$.

The function $s_{n}$, as defined above, admits partial derivatives of all
orders that are uniformly bounded. In particular, denoting for any
$\mathbf{u}\in\{(i,j)\}_{1\leq j\leq i\leq n}$, $\partial_{\mathbf{u}}s_{n}$
for $\partial s_{n}/\partial x_{\mathbf{u}}$, the following upper bounds hold:
for any $\mathbf{u},\mathbf{v},\mathbf{w}$ in $\{(i,j)\}_{1\leq j\leq i\leq
n}$, there exist universal positive constants $c_{1},c_{2}$ and $c_{3}$
depending only on the imaginary part of $z$ such that
\begin{equation}
|\partial_{\mathbf{u}}s_{n}|\leq\frac{c_{1}}{n^{3/2}}\,,\,|\partial
_{\mathbf{u}}\partial_{\mathbf{v}}s_{n}|\leq\frac{c_{2}}{n^{2}}\,\text{ and
}\,|\partial_{\mathbf{u}}\partial_{\mathbf{v}}\partial_{\mathbf{w}}s_{n}%
|\leq\frac{c_{3}}{n^{5/2}}\,. \label{boundsd}%
\end{equation}
(See Chatterjee (2006)). In addition, concerning the partial derivatives of
second order, the following lemma will be also useful.

\begin{lemma}
\label{derivatives} Let $z\in\mathbb{C}^{+}$ and $s_{n}:=s_{n,z}$ be defined
by \eqref{deffA}. Let $(a_{ij})_{1\leq j\leq i\leq n}$ and $(b_{ij})_{1\leq
j\leq i\leq n}$ be real numbers. Then, there exists an universal positive
constant $c_{4}$ depending only on the imaginary part of $z$ such that for any
subset ${\mathcal{I}}_{n}$ of $\{(i,j)\}_{1\leq j\leq i\leq n}$ and any
element ${\mathbf{x}}$ of ${\mathbb{R}}^{N}$,
\[
\Big |\sum_{\mathbf{u}\in{\mathcal{I}}_{n}}\sum_{\mathbf{v}\in{\mathcal{I}%
}_{n}}a_{\mathbf{u}}b_{\mathbf{v}}\partial_{\mathbf{u}}\partial_{\mathbf{v}%
}s_{n}({\mathbf{x}})\Big |\leq\frac{c_{4}}{n^{2}}\Big (\sum_{\mathbf{u}%
\in{\mathcal{I}}_{n}}a_{\mathbf{u}}^{2}\sum_{\mathbf{v}\in{\mathcal{I}}_{n}%
}b_{\mathbf{v}}^{2}\Big )^{1/2}\,.
\]

\end{lemma}

\noindent\textbf{Proof.} Setting $G=(A(\mathbf{x})-z{\mathbf{{I}}_{n}})^{-1}$,
we have
\[
\partial_{\mathbf{u}}\partial_{\mathbf{v}}s_{n}=\frac{1}{n}\mathrm{Tr}%
(G\partial_{\mathbf{u}}AG\partial_{\mathbf{v}}AG)+\frac{1}{n}\mathrm{Tr}%
(G\partial_{\mathbf{v}}AG\partial_{\mathbf{u}}AG)\,.
\]
(See the equality (20) in Chatterjee (2006)). Whence, with the notations
\[
{\widetilde{A}}:=\sum_{\mathbf{u}\in{\mathcal{I}}_{n}}a_{\mathbf{u}}%
\partial_{\mathbf{u}}A\,\text{\ and }{\widetilde{B}}:=\sum_{\mathbf{u}%
\in{\mathcal{I}}_{n}}b_{\mathbf{u}}\partial_{\mathbf{u}}A\,,
\]
it follows that
\[
\sum_{\mathbf{u}\in{\mathcal{I}}_{n}}\sum_{\mathbf{v}\in{\mathcal{I}}_{n}%
}a_{\mathbf{u}}b_{\mathbf{v}}\partial_{\mathbf{u}}\partial_{\mathbf{v}}%
s_{n}=\frac{2}{n}\mathrm{Tr}(G^{2}{\widetilde{A}}G{\widetilde{B}})\,.
\]
Recall now the following facts: Let $B$ and $C$ be two complex valued matrices
of order $n$. Then, $|\mathrm{Tr}(BC)|\leq\Vert B\Vert_{2}\Vert C\Vert_{2}$
where $\Vert B\Vert_{2}^{2}=\sum_{i=1}^{n}\sum_{j=1}^{n}b_{ij}^{2}$ (the
$b_{ij}$'s being the entries of $B$) and $\max\{\Vert BC\Vert_{2},\Vert
CB\Vert_{2}\}\leq\max_{1\leq i\leq n}|\eta_{i}|.\Vert C\Vert_{2}$ if $B$
admits a spectral decomposition with eigenvalues $\eta_{1},\dots,\eta_{n}$.
Therefore using the above facts, together with the facts that $(\partial
_{\mathbf{u}}A)_{ij}=n^{-1/2}$ if $(i,j)=\mathbf{u}$ or $(j,i)=\mathbf{u}$ and
$0$ otherwise, and that $G$ admits a spectral decomposition with eigenvalues
bounded by $1/y$ with $y=\mathrm{Im}(z)$, we get
\[
\frac{1}{n}|\mathrm{Tr}(G^{2}{\widetilde{A}}G{\widetilde{B}})|\leq\Vert
G^{2}\widetilde{A}\Vert_{2}\Vert G\widetilde{B}\Vert_{2}\leq\frac{1}{y^{2}%
}\frac{2}{n^{2}}\Big (\sum_{\mathbf{u}\in{\mathcal{I}}_{n}}a_{\mathbf{u}}%
^{2}\sum_{\mathbf{v}\in{\mathcal{I}}_{n}}b_{\mathbf{v}}^{2}\Big )^{1/2}\,,
\]
proving the lemma. \ $\ \ \ \lozenge$

\medskip

Another key result we use for dealing with Gaussian vectors is:

\begin{lemma}
\label{proj gaussian}Let $X=(X_{k})_{1\leq k\leq n}$ and $Y=(Y_{k})_{1\leq
k\leq n}$ be two vectors in ${\mathbb{L}}^{2}$ which have the same covariance
structure. Assume in addition that $Y$ is Gaussian. Then, for all $u\leq k$ we
have%
\[
\Vert{\mathbb{E}}(Y_{k}|\mathcal{F}_{u}^{Y})\Vert_{2}\leq\Vert{\mathbb{E}%
}(X_{k}|\mathcal{F}_{u}^{X})\Vert_{2} \, ,
\]
where $\mathcal{F}_{u}^{Y}=\sigma(Y_{i},i\leq u)$ and $\mathcal{F}_{u}%
^{X}=\sigma(X_{i},i\leq u)$.
\end{lemma}

\noindent\textbf{Proof.} To prove the inequality above, it suffices to notice
the following facts. Let
\[
{\mathcal{V}}_{u}^{Y}=\overline{\mathrm{span}}(1,(Y_{j}\,,\,1\leq j\leq
u))\,\text{ and }\,{\mathcal{V}}_{u}^{X}=\overline{\mathrm{span}}%
(1,(X_{j}\,,\,1\leq j\leq u)) \, ,
\]
where the closure is taken in ${\mathbb{L}}^{2}$. Denote by $\Pi
_{{\mathcal{V}}_{u}^{Y}}(\cdot)$ the orthogonal projection on ${\mathcal{V}%
}_{u}^{Y}$ and by $\Pi_{{\mathcal{V}}_{u}^{X}}(\cdot)$ the orthogonal
projection on ${\mathcal{V}}_{u}^{X}$. Since $(Y_{j})_{1\leq j\leq n}$ is a
Gaussian vector ${\mathbb{E}}(Y_{k}|{\mathcal{F}}_{u}^{Y})=\Pi_{{\mathcal{V}%
}_{u}^{Y}}(Y_{k})\,$ a.s. and in ${\mathbb{L}}^{2}$. On another hand, since
$(Y_{k})_{1\leq k\leq n}$ has the same covariance structure as $(X_{k})_{1\leq
k\leq n}$, we observe that
\[
\Vert\Pi_{{\mathcal{V}}_{u}^{Y}}(Y_{k})\Vert_{2}=\Vert\Pi_{{\mathcal{V}%
_{u}^{X}}^{{}}}(X_{k})\Vert_{2}\,.
\]
But, by the definition of the conditional expectation, $\Vert X_{k}%
-{\mathbb{E}}(X_{k}|{\mathcal{F}}_{u}^{X})\Vert_{2}\leq\Vert X_{k}%
-\Pi_{{\mathcal{V}}_{u}^{X}}(X_{k})\Vert_{2}$. Hence, by Pythagora's theorem,
\[
\Vert\Pi_{{\mathcal{V}}_{u}^{X}}(X_{k})\Vert_{2}\leq\Vert{\mathbb{E}}%
(X_{k}|\mathcal{F}_{u}^{X})\Vert_{2}\,.
\]
Combining all the observations above, the lemma follows. \ $\ \ \ \lozenge$

\medskip

Our next proposition gives in particular a well-known linear representation
for stationary Gaussian processes which have a spectral density. It can be
found in Varadhan (Ch 6, Section 6.6., (2001)); see also Fact 3.1 in
Chakrabarty \textit{et al.} (2014).

\begin{proposition}
\label{prop GaussLinear} Let $f$ be the spectral density on $[-\pi,\pi]$ of a
real-valued ${\mathbb{L}}^{2}$-stationary process. For any $k\in{\mathbb{Z}}$,
let
\begin{equation}
a_{k}=\frac{1}{\sqrt{2\pi}}\int_{-\pi}^{\pi}e^{{\mathrm{i}}kx}\sqrt{f(x)}dx\,.
\label{defakf}%
\end{equation}
Then the $a_{k}$'s are real numbers satisfying $\sum_{k\in{\mathbb{Z}}}%
a_{k}^{2}<\infty$. In addition, if we define for any $k\in{\mathbb{Z}}$,
\[
Y_{k}=\sum_{j\in{\mathbb{Z}}}a_{j}\xi_{k-j}\,,
\]
where $(\xi_{j})_{j\in{\mathbb{Z}}}$ is a sequence of i.i.d. standard normal
real-valued random variables, then $(Y_{k})_{k\in{\mathbb{Z}}}$ is a centered
real-valued stationary Gaussian process with spectral density $f$ on
$[-\pi,\pi]$.
\end{proposition}

\subsection{Proof of Theorem \ref{mainindependentrows}}

The proof of this theorem requires several steps. First we reduce the problem
to studying the difference of expected values. Next, in order to weaken the
dependence, we partition the variables in each row in big and small blocks.
The big blocks are approximated by vector valued martingale differences. Then,
we replace one by one these martingale differences by Gaussian random vectors
having the same covariance structure with the help of a blockwise
Lindeberg-type method.

All along the proof $z=x+{\mathrm{i}}y$ will be a complex number in
${\mathbb{C}}^{+}$. \smallskip\ Also, the notation $a\ll b$ means that there
is a constant $C$ depending only on $\operatorname{Im}z=y$ such that $a\leq
Cb.$

\bigskip

\textbf{Step 1: Reduction of the problem to a difference of expected values.}

Since the random vectors $(R_{i})_{1\leq i\leq n},$ where $R_{i}%
=(X_{ij})_{1\leq j\leq i},$ are mutually independent, it is well-known (see
for instance the arguments in the proof on page 34 in Bai-Silverstein, 2010)
that
\begin{equation}
\text{$S^{\mathbb{X}_{n}}(z)-\mathbb{E}S^{\mathbb{X}_{n}}(z)\rightarrow0$
a.s.} \label{fact1TH2}%
\end{equation}
Hence, in order to prove Theorem \ref{mainindependentrows}, it suffices to
show that
\begin{equation}
\mathbb{E}S^{\mathbb{X}_{n}}(z)-\mathbb{E}S^{\mathbb{Y}_{n}}(z)\rightarrow0.
\label{fact3TH2}%
\end{equation}
To prove the above convergence, notice that there is no loss of generality in
assuming that the entries $(Y_{\mathbf{u}})$ of $\mathbf{Y}_{n}$ are
independent of the entries $(X_{\mathbf{u}})$ of $\mathbf{X}_{n}$. Therefore,
from now on, we assume that $\mathbf{Y}_{n}$ is a symmetric matrix constructed
from a real-valued centered Gaussian random field $(Y_{{\mathbf{u}}})$ having
the same covariance structure as $(X_{\mathbf{u}})$ and independent of
$(X_{\mathbf{u}})$.

We write $S^{\mathbb{X}_{n}}(z)$ and $S^{\mathbb{Y}_{n}}(z)$ as a function of
the entries on and below the diagonal, arranged row after row. More exactly,
using the notation \eqref{deffA}, we write
\[
\mathbb{\ }S^{\mathbb{X}_{n}}(z)=s_{n}\big (L^{X}\big )\text{ and }%
\mathbb{\ }S^{\mathbb{Y}_{n}}(z)=s_{n}\big (L^{Y}\big )\,,
\]
where $L^{X}=(L_{i}^{X})_{1\leq i\leq n}$ and $L^{Y}=(L_{i}^{Y})_{1\leq i\leq
n}$ with $L_{i}^{X}=(X_{i1},\dots,X_{ii})$ and $L_{i}^{Y}=(Y_{i1},\dots
,Y_{ii})$. Also, in the sequel, to further simplify the notation we shall skip
the index $n$ from $s_{n}$ and we put $s=$ $s_{n}:=s_{n,z}.$

\medskip

\noindent\textbf{Step 2: Martingale approximation.}

We shall introduce a martingale structure on each row. We start from the
celebrated Bernstein big and small blocks argument which weakens the
dependence. We partition the variables in each row in big and small blocks and
show that the variables in large blocks have a dominant contribution. The
large blocks are then decomposed in martingale differences and a rest which
also has a smaller contribution.

Let $p$ and $q$ be two integers fixed for the moment. Fix $i$ in
$\{1,\dots,n\}$ and let $k_{i}=[i/(p+q)]$. We partition the set $\{1,\dots
,i\}$ in big and small blocks with the following restriction: a big block of
size $p$ is followed by a small block of size $q$. We shall have the set of
indexes $I_{i1},J_{i1},I_{i2},J_{i2},....,I_{ik_{i}},J_{i,k_{i}+1}$ where each
index set $I_{ij}$ is of size $p$, each index set $J_{ij}$ is of size $q$ and
the last block has a size at most $p+q$. More precisely, for any $i$ in
$\{1,\dots,n\}$ and for any $j\in\{1,\dots,k_{i}\}$,
\[
I_{j}^{\prime}=\{(j-1)(p+q)+1,\dots,(j-1)(p+q)+p\}\,\text{ and }%
\,J_{j}^{\prime}=\{(j-1)(p+q)+p+1,\dots,j(p+q)\}\,.
\]
and%
\[
I_{ij}=\{(i,k);k\in I_{j}^{\prime}\,\}\text{ and }\,J_{ij}=\{(i,k);k\in
J_{j}^{\prime}\,\}\,.
\]
Corresponding to this index decomposition, the vectors $L_{j}^{X}$ and
$L_{j}^{Y}$ are partitioned in $k_{i}+1$ consecutive vectors. Setting
\[
B_{ij}=(X_{{\mathbf{u}}})_{{\mathbf{u}}\in I_{ij}}\,,\,b_{ij}=(X_{{\mathbf{u}%
}})_{{\mathbf{u}}\in J_{ij}}\,,\,B_{ij}^{\ast}=(Y_{{\mathbf{u}}}%
)_{{\mathbf{u}}\in I_{ij}}\,\text{ and }\,b_{ij}^{\ast}=(Y_{{\mathbf{u}}%
})_{{\mathbf{u}}\in J_{ij}}%
\]
we write
\[
L_{i}^{X}=(B_{i1},b_{i1},B_{i2},b_{i2},\cdots,B_{ik_{i}},b_{ik_{i}}%
,b_{i,k_{i}+1})\,\text{ and }\,L_{i}^{Y}=(B_{i1}^{\ast},b_{i1}^{\ast}%
,B_{i2}^{\ast},b_{i2}^{\ast},\cdots,B_{ik_{i}}^{\ast},b_{ik_{i}}^{\ast
},b_{i,k_{i}+1}^{\ast})\,.
\]
We introduce now the following vectors
\[
B_{i}^{X}=(B_{i1},{\mathbf{0}}_{q},B_{i2},{\mathbf{0}}_{q},\cdots,B_{ik_{i}%
},{\mathbf{0}}_{q},{\mathbf{0}}_{r})\,\text{ and }\,B_{i}^{Y}=(B_{i1}^{\ast
},{\mathbf{0}}_{q},B_{i2}^{\ast},{\mathbf{0}}_{q},\cdots,B_{ik_{i}}^{\ast
},{\mathbf{0}}_{q},{\mathbf{0}}_{r})\,,
\]
where $r=i-k_{i}(p+q)$. Note that $B_{i}^{X}$ (resp. $B_{i}^{Y}$) is derived
from $L_{i}^{X}$ (resp. $L_{i}^{Y}$) where we replace the variables in
$b_{ij}$ (resp. $b_{ij}^{\ast}$) by $0$'s. In addition, for $A$ a positive
real, fixed for the moment, we set for any $\mathbf{u}\in{\mathbb{Z}}^{2}$
\[
{\widetilde{X}}_{\mathbf{u}}:=X_{\mathbf{u}}I(|X_{\mathbf{u}}|\leq A)\,,
\]
and, for any $i\in\{1,\dots,n\}$,
\[
{\widetilde{B}}_{i}^{X}=({\widetilde{B}}_{i1},{\mathbf{0}}_{q},{\widetilde{B}%
}_{i2},{\mathbf{0}}_{q},\cdots,{\widetilde{B}}_{ik_{i}},{\mathbf{0}}%
_{q},{\mathbf{0}}_{r})\,\text{ where }\,{\widetilde{B}}_{ij}=({\widetilde{X}%
}_{\mathbf{u}})_{\mathbf{u}\in I_{ij}}\text{ for $j\in\{1,\dots,k_{i}\}$.}%
\]
Next, for any $i\in\{1,\dots,n\}$, we consider the sigma algebras
$\mathcal{F}_{i0}^{X}=\mathcal{F}_{i0}^{Y}=\{\emptyset,\Omega\}$ and for
$1\leq\ell\leq k_{i}$, $\mathcal{F}_{i\ell}^{X}=\sigma(B_{ij};1\leq j\leq
\ell)$ and $\mathcal{F}_{i\ell}^{Y}=\sigma(B_{ij}^{\ast};1\leq j\leq\ell)$.
Then, for any $\ell\in\{1,\dots,k_{i}\}$, we define
\begin{equation}
{\widetilde{D}}_{i\ell}={\widetilde{B}}_{i\ell}-\mathbb{E}({\widetilde{B}%
}_{i\ell}|\mathcal{F}_{i,\ell-1}^{X})\,, \label{defD}%
\end{equation}
and
\begin{equation}
{D}_{i\ell}^{\ast}={B}_{i\ell}^{\ast}-\mathbb{E}({B}_{i\ell}^{\ast
}|\mathcal{F}_{i,\ell-1}^{Y})\,. \label{defDG}%
\end{equation}
By $\mathbb{E}({\widetilde{B}}_{i\ell}|\mathcal{F}_{i,\ell-1}^{X})$ (resp.
$\mathbb{E}({B}_{i\ell}^{\ast}|\mathcal{F}_{i,\ell-1}^{Y})$) we understand a
vector of dimension $p$ where each component is a component of the vector
${\widetilde{B}}_{i\ell}$ (resp. ${B}_{i\ell}^{\ast}$) conditioned with
respect to $\mathcal{F}_{i,\ell-1}^{X}$ (resp. $\mathcal{F}_{i,\ell-1}^{Y}$).
Note that $({\widetilde{D}}_{i\ell})_{1\leq\ell\leq k_{i}}$ and $({D}_{i\ell
}^{\ast})_{1\leq\ell\leq k_{i}}$ are vector valued martingale differences
adapted respectively to $(\mathcal{F}_{i\ell}^{X})_{1\leq\ell\leq k_{i}}$ and
$(\mathcal{F}_{i\ell}^{Y})_{1\leq\ell\leq k_{i}}$. We then define the vectors
${\widetilde{D}}_{i}^{X}$ and $D_{i}^{Y}$ with dimension $i$ and with a
similar structure as $B_{i}^{X}$ as follows:
\begin{equation}
{\widetilde{D}}_{i}^{X}=({\widetilde{D}}_{i1},{\mathbf{0}}_{q},{\widetilde{D}%
}_{i2},{\mathbf{0}}_{q},\cdots,{\widetilde{D}}_{ik_{i}},{\mathbf{0}}%
_{q},{\mathbf{0}}_{r})\,\text{ and }\,{D}_{i}^{Y}=({D}_{i1}^{\ast}%
,{\mathbf{0}}_{q},{D}_{i2}^{\ast},{\mathbf{0}}_{q},\cdots,{D}_{ik_{i}}^{\ast
},{\mathbf{0}}_{q},{\mathbf{0}}_{r})\,. \label{linesD}%
\end{equation}
Setting ${\widetilde{D}}^{X}=({\widetilde{D}}_{i}^{X})_{1\leq i\leq n}$, we
first compare ${\mathbb{E}}s(L^{X})$ to ${\mathbb{E}}s({\widetilde{D}}^{X})$.
\ We write%
\[
{\mathbb{E}}s(L^{X})-{\mathbb{E}}s({\widetilde{D}}^{X})=\mathbb{E}\Delta
_{1}(s)+\mathbb{E}\Delta_{2}(s)+\mathbb{E}\Delta_{3}(s)\,,
\]
where
\[
\Delta_{1}(s)=s(L^{X})-s(B^{X})\,,\,\Delta_{2}(s)=s(B^{X})-s({\widetilde{B}%
}^{X})
\]
and
\[
\Delta_{3}(s)=s({\widetilde{B}}^{X})-s({\widetilde{D}}^{X})\,,
\]
with the notations $B^{X}=(B_{i}^{X})_{1\leq i\leq n}$ and ${\widetilde{B}%
}^{X}=({\widetilde{B}}_{i}^{X})_{1\leq i\leq n}$. To control each of the
$\mathbb{E}\Delta_{i}(s)$ for $i=1,2,3$, we apply Lemma \ref{lmagotze}.
Therefore, we get
\[
|\mathbb{E}\Delta_{1}(s)|^{2}\ll\sum_{i=1}^{n}\sum_{j=1}^{k_{i}+1}%
\sum_{\mathbf{u}\in J_{ij}}{\mathbb{E}}(X_{\mathbf{u}}^{2})\ll\Big (\frac
{q}{p}+\frac{q+p}{n}\Big )\sigma^{2}\,,
\]%
\[
|\mathbb{E}\Delta_{2}(s)|^{2}\ll\frac{1}{n^{2}}\sum_{i=1}^{n}\sum_{j=1}%
^{k_{i}}\sum_{\mathbf{u}\in I_{ij}}{\mathbb{E}}(X_{\mathbf{u}}^{2}%
I(|X_{\mathbf{u}}|>A))\ll L(A)\,,
\]
and
\begin{align*}
|\mathbb{E}\Delta_{3}(s)|^{2}  &  \ll\frac{1}{n^{2}}\sum_{i=1}^{n}\sum
_{j=1}^{k_{i}}\sum_{\mathbf{u}\in I_{ij}}\Vert{\mathbb{E}}({\widetilde{X}%
}_{\mathbf{u}}|\mathcal{F}_{i,j-1}^{X})\Vert_{2}^{2}\leq2\left(
L(A)+\max_{1\leq j\leq i\leq n}\Vert{\mathbb{E}}({X}_{ij}|\mathcal{G}%
_{i,j-q})\Vert_{2}^{2}\right) \\
&  \ll\left(  L(A)+\eta_{q}^{2}\right)  \,.
\end{align*}
We proceed in a similar way with the matrix $\mathbb{Y}_{n}$. Therefore,
setting ${D}^{Y}=({D}_{i}^{Y})_{1\leq i\leq n}$, we write%
\[
\mathbb{E}s(L^{Y})-\mathbb{E}s(D^{Y})=\mathbb{E}\Delta_{1}^{\prime
}(s)+\mathbb{E}\Delta_{2}^{\prime}(s),
\]
with the notations
\[
\Delta_{1}^{\prime}(s)=s(L^{Y})-s(B^{Y})\,\text{ and }\,\Delta_{2}^{\prime
}(s)=s(B^{Y})-s({D}^{Y}),
\]
where $B^{Y}=(B_{i}^{Y})_{1\leq i\leq n}$. Applying Lemma \ref{lmagotze} and
using the fact that $(Y_{\mathbf{u}})$ has the same covariance structure as
$(X_{\mathbf{u}})$, we derive
\[
|\mathbb{E}\Delta_{1}^{\prime}(s_{n})|^{2}\ll\Big (\frac{q}{p}+\frac{q+p}%
{n}\Big )\sup{\mathbb{E}}(Y_{\mathbf{u}}^{2})\ll\Big (\frac{q}{p}+\frac
{q+p}{n}\Big )\sigma^{2}\,.
\]
On another hand, Lemmas \ref{lmagotze} and \ref{proj gaussian} imply that
\begin{align*}
|\mathbb{E}\Delta_{2}^{\prime}(s)|^{2}  &  \ll\frac{1}{n^{2}}\sum_{i=1}%
^{n}\sum_{j=1}^{k_{i}}\sum_{\mathbf{u}\in I_{ij}}\Vert{\mathbb{E}}%
({Y}_{\mathbf{u}}|\mathcal{F}_{i,j-1}^{Y})\Vert_{2}^{2}\ll\frac{1}{n^{2}}%
\sum_{i=1}^{n}\sum_{j=1}^{k_{i}}\sum_{\mathbf{u}\in I_{ij}}\Vert{\mathbb{E}%
}({X}_{\mathbf{u}}|\mathcal{F}_{i,j-1}^{X})\Vert_{2}^{2}\\
&  \ll\max_{1\leq j\leq i\leq n}\Vert{\mathbb{E}}({X}_{ij}|\mathcal{G}%
_{i,j-q})\Vert_{2}^{2}\,\ll\eta_{q}^{2}.
\end{align*}
Overall we have the decomposition%
\begin{equation}
\mathbb{E}S^{\mathbb{X}_{n}}(z)-\mathbb{E}S^{\mathbb{Y}_{n}}(z)=\mathbb{E}%
s({\widetilde{D}}^{X})-\mathbb{E}s(D^{Y})+E_{n}(p,q,A)\,,
\label{mart decompose}%
\end{equation}
with
\[
|E_{n}(p,q,A)|\ll\Big (\Big (\frac{q}{p}+\frac{q+p}{n}\Big )^{1/2}%
\sigma+L^{1/2}(A)+\eta_{q}\Big )\,.
\]

\medskip

\noindent\textbf{Step 3: The study of }$\mathbb{E}s({\widetilde{D}}%
^{X})-\mathbb{E}s(D^{Y})$.

To study $\mathbb{E}s({\widetilde{D}}^{X})-\mathbb{E}s(D^{Y})$ we first
decompose the difference according to the rows and after that we study the
rows separately. With this aim we introduce a telescoping sum where each term
is a difference of two functions whose arguments differ only by one row.
Namely we write
\[
\mathbb{E}s({\widetilde{D}}^{X})-\mathbb{E}s(D^{Y})=\sum_{i=1}^{n}\left(
\mathbb{E}s\big ({\widetilde{D}}_{[1,i-1]}^{X},{\widetilde{D}}_{i}^{X}%
,{D}_{[i+1,n]}^{Y}\big )-\mathbb{E}s\big ({\widetilde{D}}_{[1,i-1]}^{X}%
,{D}_{i}^{Y},{D}_{[i+1,n]}^{Y}\big )\right)
\]
where ${\widetilde{D}}_{[a,b]}^{X}=({\widetilde{D}}_{a}^{X},...{\widetilde{D}%
}_{b}^{X})$ and ${D}_{[a,b]}^{Y}=({D}_{a}^{Y},...{D}_{b}^{Y})$ with
${\widetilde{D}}_{i}^{X}$ and ${D}_{i}^{Y}$ defined in \eqref{linesD}. Now for
every $i$ fixed denote by
\[
s_{i}({\mathbf{x}}):=s\big ({\widetilde{D}}_{[1,i-1]}^{X},{\mathbf{x}}%
,{D}_{[i+1,n]}^{Y}\big )\,\,.
\]
Note that $s_{i}$ is a random function from ${\mathbb{R}}^{i}$ to
${\mathbb{C}}$. With this notation%
\[
\mathbb{E}s({\widetilde{D}}^{X})-\mathbb{E}s(D^{Y})=\sum_{i=1}^{n}%
\mathbb{E}(s_{i}({\widetilde{D}}_{i}^{X})-s_{i}(D_{i}^{Y}))\,.
\]
From now on, for easier notation, it will be convenient to extend the vectors
$({\widetilde{D}}_{i\ell})_{1\leq\ell\leq k_{i}}$ and $({D}_{i\ell}^{\ast
})_{1\leq\ell\leq k_{i}}$ defined in \eqref{defD} and \eqref{defDG} as
follows:
\begin{equation}
{\widetilde{D}}_{i\ell}^{\prime}=({\widetilde{D}}_{i\ell},{\mathbf{0}}%
_{q})\,\text{ and }\,{D^{\prime}}_{i\ell}^{\ast}=({\widetilde{D}}_{i\ell
}^{\ast},{\mathbf{0}}_{q})\,\text{ for $1\leq\ell\leq k_{i}-1$}
\label{defDext1}%
\end{equation}
and
\begin{equation}
{\widetilde{D}}_{ik_{i}}^{\prime}=({\widetilde{D}}_{ik_{i}},{\mathbf{0}}%
_{q+r})\,\text{ and }\,{D^{\prime}}_{ik_{i}}^{\ast}=({D}_{ik_{i}}^{\ast
},{\mathbf{0}}_{q+r})\,. \label{defDext2}%
\end{equation}
With these notations, as in the Lindeberg's method, we write now another
telescoping sum where we change one by one the vectors ${\widetilde{D}}%
_{i\ell}^{\prime}$ by ${D^{\prime}}_{i\ell}^{\ast}$ in the argument of $s_{i}%
$. With this aim we write
\begin{align}
s_{i}({\widetilde{D}}_{i}^{X})  &  -s_{i}(D_{i}^{Y})=s_{i}({\widetilde{D}%
}_{i1}^{\prime},\dots,{\widetilde{D}}_{ik_{i}}^{\prime})-s_{i}({D^{\prime}%
}_{i1}^{\ast},\dots,{D^{\prime}}_{ik_{i}}^{\ast})\nonumber\\
&  =\sum_{u=1}^{k_{i}}\big (s_{i}({\widetilde{D}}_{i,[1,u-1]}^{\prime
},{\widetilde{D}}_{iu}^{\prime},{D^{\prime}}_{i,[u+1,k_{i}]}^{\ast}%
)-s_{i}({\widetilde{D}}_{i,[1,u-1]}^{\prime},{D^{\prime}}_{iu}^{\ast
},{D^{\prime}}_{i,[u+1,k_{i}]}^{\ast}\big )\nonumber\\
&  :=\sum_{u=1}^{k_{i}}\big (s_{i,u}({\widetilde{D}}_{iu}^{\prime}%
)-s_{i,u}({D^{\prime}}_{iu}^{\ast})\big )\,, \label{TErows}%
\end{align}
where ${\widetilde{D}}_{i,[k,\ell]}^{\prime}:=({\widetilde{D}}_{ik}^{\prime
},\dots,{\widetilde{D}}_{i\ell}^{\prime})$ and ${D^{\prime}}_{i,[k,\ell
]}^{\ast}:=({D^{\prime}}_{ik}^{\ast},\dots,{D^{\prime}}_{i\ell}^{\ast})$. Note
that the $s_{iu}$'s defined above are random functions from ${\mathbb{R}%
}^{p+q}$ to ${\mathbb{C}}$ if $1\leq u\leq k_{i}-1$ and from ${\mathbb{R}%
}^{q+r}$ to ${\mathbb{C}}$ if $u=k_{i}$ (where $r=i-k_{i}(p+q)$).

We shall treat separately each term in the sum (\ref{TErows}) corresponding to
the $i$-th row. So, in the following, $i$ is fixed. To facilitate the study of
this difference we introduce some auxiliary terms:%
\[
s_{iu}({\widetilde{D}}_{iu}^{\prime})-s_{iu}({D^{\prime}}_{iu}^{\ast}%
)=s_{iu}({\widetilde{D}}_{iu}^{\prime})-s_{iu}({\mathbf{0}})+s_{iu}%
({\mathbf{0}})-s_{iu}({D^{\prime}}_{iu}^{\ast})\,.
\]
Denote by $d_{iu}^{(j)}$ the $j$-th component of the vector ${\widetilde{D}%
}_{iu}^{\prime}$. Using Taylor's expansion of order three, we get
\begin{equation}
s_{iu}({\widetilde{D}}_{iu}^{\prime})-s_{iu}({\mathbf{0}})={\widetilde{R}}%
_{1}+{\widetilde{R}}_{2}+{\widetilde{R}}_{3}\,, \label{TErows1}%
\end{equation}
where
\[
{\widetilde{R}}_{1}=\sum_{j=1}^{p}d_{iu}^{(j)}\partial_{j}s_{iu}({\mathbf{0}%
})\,,\text{ }\ \text{ }\,{\widetilde{R}}_{2}=\frac{1}{2}\Big (\sum_{j=1}%
^{p}d_{iu}^{(j)}\partial_{j}\Big )^{2}s_{iu}({\mathbf{0}})
\]
and
\[
{\widetilde{R}}_{3}=\frac{1}{6}\Big (\sum_{j=1}^{p}d_{iu}^{(j)}\partial
_{j}\Big )^{3}s_{iu}(\theta{\widetilde{D}}_{iu}^{\prime})\,\text{ with
$\theta\in]0,1[$}\,.
\]
Similarly, if we denote by $g_{iu}^{(j)}$ the $j$-th component of the vector
${D^{\prime}}_{iu}^{\ast}$, we get
\begin{equation}
s_{iu}({D^{\prime}}_{iu}^{\ast})-s_{iu}({\mathbf{0}})={R}_{1}^{\ast}+{R}%
_{2}^{\ast}+{R}_{3}^{\ast}\,, \label{TErows2}%
\end{equation}
where
\[
R_{1}^{\ast}=\sum_{j=1}^{p}g_{iu}^{(j)}\partial_{j}s_{iu}({\mathbf{0}%
})\,\text{ and }\,R_{2}^{\ast}=\frac{1}{2}\Big (\sum_{j=1}^{p}g_{iu}%
^{(j)}\partial_{j}\Big )^{2}s_{iu}({\mathbf{0}})
\]
and
\[
R_{3}^{\ast}=\frac{1}{6}\Big (\sum_{j=1}^{p}g_{iu}^{(j)}\partial_{j}%
\Big )^{3}s_{iu}(\theta{D^{\prime}}_{iu}^{\ast})\,\text{ with $\theta\in
]0,1[$}\,.
\]
Now notice that, for any $u\in\{1,\dots,k_{i}\}$ and any $j\in\{1,\dots,p\}$,
\begin{equation}
d_{iu}^{(j)}={\widetilde{X}}_{i,(u-1)(p+q)+j}-{\mathbb{E}}({\widetilde{X}%
}_{i,(u-1)(p+q)+j}|{\mathcal{F}}_{i,u-1}^{X}):={\widetilde{X}}_{iu}%
^{(j)}-{\mathbb{E}}({\widetilde{X}}_{iu}^{(j)}|{\mathcal{F}}_{i,u-1}^{X})\,,
\label{componentd}%
\end{equation}
and
\begin{equation}
g_{iu}^{(j)}={Y}_{i,(u-1)(p+q)+j}-{\mathbb{E}}(Y_{i,(u-1)(p+q)+j}%
|{\mathcal{F}}_{i,u-1}^{Y}):=Y_{iu}^{(j)}-{\mathbb{E}}(Y_{iu}^{(j)}%
|{\mathcal{F}}_{i,u-1}^{Y})\,. \label{componentg}%
\end{equation}
Therefore
\[
\Vert d_{iu}^{(j)}\Vert_{3}^{3}\leq2^{3}\Vert{\widetilde{X}}_{iu}^{(j)}%
\Vert_{3}^{3}\ll A\sigma^{2}\,,
\]
and since $\mathbf{Y}_{n}$ has the same covariance structure as $\mathbf{X}%
_{n}$ and is a Gaussian vector,
\[
\Vert g_{iu}^{(j)}\Vert_{3}^{3}\leq2^{3}\Vert Y_{iu}^{(j)}\Vert_{3}^{3}%
\leq2^{4}\Vert Y_{iu}^{(j)}\Vert_{2}^{3}\ll\sigma^{3}.
\]
Taking into account the two previous inequalities and the upper bound on the
partial derivatives of order three of $s$ given in \eqref{boundsd}, we infer
that
\begin{equation}
|{\mathbb{E}}({\widetilde{R}}_{3})+{\mathbb{E}}(R_{3}^{\ast})|\ll\frac
{1}{n^{5/2}}p^{3}\sigma^{2}(A+\sigma)\,. \label{UBR3}%
\end{equation}
On another hand, we notice that for any $j,\ell$ in $\{1,\dots,p\}$,
$\partial_{j}s_{iu}({\mathbf{0}})$ and $\partial_{j}\partial_{\ell}%
s_{iu}({\mathbf{0}})$ are complex-valued random variables measurable with
respect to the sigma algebra ${\mathcal{H}}_{i,u}$ defined by
\begin{equation}
{\mathcal{H}}_{i,u}={\mathcal{F}}_{i,u-1}^{X}\vee\sigma\big ((L_{j}%
^{X})_{1\leq j\leq i-1},(L_{k}^{Y})_{i+1\leq k\leq n}\big )\vee\sigma
\big ({D}_{i,u+1}^{\ast},\dots,{D}_{ik_{i}}^{\ast}\big )\,.
\label{definitionHX}%
\end{equation}
Hence
\[
{\mathbb{E}}({\widetilde{R}}_{1})=\sum_{j=1}^{p}{\mathbb{E}}\big (\partial
_{j}s_{iu}({\mathbf{0}}){\mathbb{E}}(d_{iu}^{(j)}|{\mathcal{H}}_{i,u}%
)\big )\,,
\]
and
\[
{\mathbb{E}}({\widetilde{R}}_{2})=\frac{1}{2}\sum_{j,\ell=1}^{p}{\mathbb{E}%
}\big (\partial_{j}\partial_{\ell}s_{iu}({\mathbf{0}}){\mathbb{E}}%
(d_{iu}^{(j)}d_{iu}^{(\ell)}|{\mathcal{H}}_{i,u})\big )\,.
\]
Since the rows of ${\mathbf{X}}_{n}$ are assumed to be independent and
${\mathbf{Y}}_{n}$ is assumed to be independent of ${\mathbf{X}}_{n}$, then
$\sigma(d_{iu}^{(1)},\dots,d_{iu}^{(p)})\vee{\mathcal{F}}_{i,u-1}^{X}$ is
independent of
\[
\sigma\big ((L_{j}^{X})_{1\leq j\leq i-1},(L_{k}^{Y})_{i+1\leq k\leq
n}\big )\vee\sigma\big ({D}_{i,u+1}^{\ast},\dots,{D}_{ik_{i}}^{\ast}\big ).
\]
Therefore, by the properties of the conditional expectation, ${\mathbb{E}%
}(d_{iu}^{(j)}|{\mathcal{H}}_{i,u})={\mathbb{E}}(d_{iu}^{(j)}|{\mathcal{F}%
}_{i,u-1}^{X})=0$ and ${\mathbb{E}}(d_{iu}^{(j)}d_{iu}^{(\ell)}|{\mathcal{H}%
}_{i,u})={\mathbb{E}}(d_{iu}^{(j)}d_{iu}^{(\ell)}|{\mathcal{F}}_{i,u-1}^{X})$.
Hence,
\begin{equation}
{\mathbb{E}}({\widetilde{R}}_{1})=0\,\text{ and }{\mathbb{E}}({\widetilde{R}%
}_{2})=\frac{1}{2}\sum_{j,\ell=1}^{p}{\mathbb{E}}\big ({\mathbb{E}}%
(d_{iu}^{(j)}d_{iu}^{(\ell)}|{\mathcal{F}}_{i,u-1}^{X})\partial_{j}%
\partial_{\ell}s_{iu}({\mathbf{0}})\big )\,. \label{R1R2}%
\end{equation}
We handle now the terms ${\mathbb{E}}(R_{1}^{\ast})$ and ${\mathbb{E}}%
(R_{2}^{\ast})$. With this aim we notice that by definition $(D_{iu}^{\ast
}:1\leq u\leq k_{i})_{1\leq i\leq n}$ is a centered Gaussian vector such that
${\mathrm{Cov}}(D_{iu}^{\ast},D_{i^{\prime}u^{\prime}}^{\ast})={\mathbf{0}%
}_{p,p}$ if $(i,u)\neq(i^{\prime},u^{\prime})$. Therefore $D_{i,u}^{\ast}$,
$i=1,\dots,n$, $u=1,\dots,k_{i}$ are centered Gaussian random variables in
${\mathbb{R}}^{p}$ which are mutually independent. In addition they are
independent of $(X_{\mathbf{u}})$. Therefore,
\begin{equation}
{\mathbb{E}}(R_{1}^{\ast})=\sum_{j=1}^{p}{\mathbb{E}}(g_{iu}^{(j)}%
){\mathbb{E}}\big (\partial_{j}s_{iu}({\mathbf{0}})\big )=0\,, \label{R1*}%
\end{equation}
and
\begin{equation}
{\mathbb{E}}(R_{2}^{\ast})=\frac{1}{2}\sum_{j,\ell=1}^{p}{\mathbb{E}}%
(g_{iu}^{(j)}g_{iu}^{(\ell)}){\mathbb{E}}\big (\partial_{j}\partial_{\ell
}s_{iu}({\mathbf{0}})\big )\,. \label{R2*}%
\end{equation}
So, starting from \eqref{TErows} and taking into account \eqref{TErows1},
\eqref{TErows2}, \eqref{UBR3}, \eqref{R1R2}, \eqref{R1*} and \eqref{R2*}, we
derive that for any $i\in\{1,\dots,n\}$,
\begin{multline}
{\mathbb{E}}\big(s_{i}({\widetilde{D}}_{i}^{X})\big )-{\mathbb{E}}%
\big(s_{i}(D_{i}^{Y})\big )\ll\Big |\sum_{u=1}^{k_{i}}\sum_{j,\ell=1}%
^{p}{\mathbb{E}}\Big (\big ({\mathbb{E}}(d_{iu}^{(j)}d_{iu}^{(\ell
)}|{\mathcal{F}}_{i,u-1}^{X})-{\mathbb{E}}(g_{iu}^{(j)}g_{iu}^{(\ell
)})\big )\partial_{j}\partial_{\ell}s_{iu}({\mathbf{0}}%
)\Big )\Big |\label{TErowsinter}\\
+\frac{1}{n^{5/2}}k_{i}p^{3}\sigma^{2}(A+\sigma)\,.
\end{multline}
We handle now the first term in the right-hand side of the above inequality.
Recalling the notations \eqref{componentd} and \eqref{componentg}, we first
write
\begin{multline*}
{\mathbb{E}}(d_{iu}^{(j)}d_{iu}^{(\ell)}|{\mathcal{F}}_{i,u-1}^{X}%
)-{\mathbb{E}}(g_{iu}^{(j)}g_{iu}^{(\ell)})={\mathbb{E}}({\widetilde{X}}%
_{iu}^{(j)}{\widetilde{X}}_{iu}^{(\ell)}|{\mathcal{F}}_{i,u-1}^{X}%
)-{\mathbb{E}}({Y}_{iu}^{(j)}{Y}_{iu}^{(\ell)})\\
-{\mathbb{E}}({\widetilde{X}}_{iu}^{(j)}|{\mathcal{F}}_{i,u-1}^{X}%
){\mathbb{E}}({\widetilde{X}}_{iu}^{(\ell)}|{\mathcal{F}}_{i,u-1}%
^{X})+{\mathbb{E}}\big ({\mathbb{E}}({Y}_{iu}^{(j)}|{\mathcal{F}}_{i,u-1}%
^{Y}){\mathbb{E}}(Y_{iu}^{(\ell)}|{\mathcal{F}}_{i,u-1}^{Y})\big )\,.
\end{multline*}
Therefore, by triangle inequality and Jensen inequality,
\begin{align}
\Big |\sum_{u=1}^{k_{i}}\sum_{j,\ell=1}^{p}{\mathbb{E}}\Big (  &
\big ({\mathbb{E}}(d_{iu}^{(j)}d_{iu}^{(\ell)}|{\mathcal{F}}_{i,u-1}%
^{X})-{\mathbb{E}}(g_{iu}^{(j)}g_{iu}^{(\ell)})\big )\partial_{j}%
\partial_{\ell}s_{iu}({\mathbf{0}}%
)\Big )\Big |\nonumber\label{TErowsinterterme1}\\
&  \leq\sum_{u=1}^{k_{i}}\sum_{j,\ell=1}^{p}\big |{\mathbb{E}}%
\big (\big ({\mathbb{E}}({\widetilde{X}}_{iu}^{(j)}{\widetilde{X}}_{iu}%
^{(\ell)}|{\mathcal{F}}_{i,u-1}^{X})-{\mathbb{E}}({Y}_{iu}^{(j)}{Y}%
_{iu}^{(\ell)})\big )\partial_{j}\partial_{\ell}s_{iu}({\mathbf{0}%
})\big )\big |\nonumber\\
&  \quad\quad+\sum_{u=1}^{k_{i}}{\mathbb{E}}\Big |\sum_{j,\ell=1}%
^{p}{\mathbb{E}}({\widetilde{X}}_{iu}^{(j)}|{\mathcal{F}}_{i,u-1}%
^{X}){\mathbb{E}}({\widetilde{X}}_{iu}^{(\ell)}|{\mathcal{F}}_{i,u-1}%
^{X})\partial_{j}\partial_{\ell}s_{iu}({\mathbf{0}})\Big |\nonumber\\
&  \quad\quad+\sum_{u=1}^{k_{i}}\Big |\sum_{j,\ell=1}^{p}{\mathbb{E}%
}\big ({\mathbb{E}}({Y}_{iu}^{(j)}|{\mathcal{F}}_{i,u-1}^{Y}){\mathbb{E}%
}(Y_{iu}^{(\ell)}|{\mathcal{F}}_{i,u-1}^{Y})\big ){\mathbb{E}}\big (\partial
_{j}\partial_{\ell}s_{iu}({\mathbf{0}})\big )\Big |\nonumber\\
&  :=T_{1}+T_{2}+T_{3}\,.
\end{align}
Let us first handle $T_{3}$. Recalling the notation \eqref{defDext1} and
\eqref{defDext2} and setting
\begin{equation}
C_{i,u}=\big ({\widetilde{D}}_{[1,i-1]}^{X},{\widetilde{D}}_{i1}^{\prime
},\dots,{\widetilde{D}}_{i,u-1}^{\prime},{\mathbf{0}},{D^{\prime}}%
_{i,u+1}^{\ast},\dots,{D^{\prime}}_{i,u_{k_{i}}}^{\ast},{D}_{[i+1,n]}%
^{Y}\big )\,, \label{defCiu}%
\end{equation}
we note that ${\mathbb{E}}({Y}_{iu}^{(j)}|{\mathcal{F}}_{i,u-1}^{Y}%
){\mathbb{E}}(Y_{iu}^{(\ell)}|{\mathcal{F}}_{i,u-1}^{Y})$ is independent of
$\partial_{j}\partial_{\ell}s_{iu}(C_{i,u})$. This is because of the
independence between ${\mathbf{Y}}_{n}$ and ${\mathbf{X}}_{n}$ together with
the independence between the vectors $({\mathbb{E}}({Y}_{iu}^{(j)}%
|{\mathcal{F}}_{i,u-1}^{Y}),{\mathbb{E}}(Y_{iu}^{(\ell)}|{\mathcal{F}}%
_{i,u-1}^{Y}))$ and $({D^{\prime}}_{i,u+1}^{\ast},\dots,{D^{\prime}%
}_{i,u_{k_{i}}}^{\ast},{D}_{[i+1,n]}^{Y})$. To prove the latter independence,
it suffices to notice that $({\mathbb{E}}({Y}_{iu}^{(j)}|{\mathcal{F}}%
_{i,u-1}^{Y}),{\mathbb{E}}(Y_{iu}^{(\ell)}|{\mathcal{F}}_{i,u-1}%
^{Y}),{D^{\prime}}_{i,u+1}^{\ast},\dots,{D^{\prime}}_{i,u_{k_{i}}}^{\ast}%
,{D}_{[i+1,n]}^{Y})$ is a Gaussian vector and that $({\mathbb{E}}({Y}%
_{iu}^{(j)}|{\mathcal{F}}_{i,u-1}^{Y}),{\mathbb{E}}(Y_{iu}^{(\ell
)}|{\mathcal{F}}_{i,u-1}^{Y}))$ and $({D^{\prime}}_{i,u+1}^{\ast}%
,\dots,{D^{\prime}}_{i,u_{k_{i}}}^{\ast},{D}_{[i+1,n]}^{Y})$ are uncorrelated.
So, we can bound $T_{3}$ as follows:
\[
T_{3}\leq\sum_{u=1}^{k_{i}}{\mathbb{E}}\Big |\sum_{j,k\in I_{u}^{\prime}}^{{}%
}{\mathbb{E}}(Y_{ij}|{\mathcal{F}}_{i,u-1}^{Y}){\mathbb{E}}(Y_{ik}%
|{\mathcal{F}}_{i,u-1}^{Y})\partial_{ij}\partial_{ik}s(C_{i,u})\Big |\,.
\]
An application of Lemma \ref{derivatives} gives
\[
\Big |\sum_{j,k\in I_{u}^{\prime}}^{{}}{\mathbb{E}}(Y_{ij}|{\mathcal{F}%
}_{i,u-1}^{Y}){\mathbb{E}}(Y_{ik}|{\mathcal{F}}_{i,u-1}^{Y})\partial
_{ij}\partial_{ik}s(C_{i,u})\Big |\ll\frac{1}{n^{2}}\sum_{j\in I_{u}^{\prime}%
}^{{}}\big ({\mathbb{E}}(Y_{ij}|{\mathcal{F}}_{i,u-1}^{Y})\big )^{2}\,.
\]
Whence, using in addition Lemma \ref{proj gaussian}, we derive
\[
T_{3}\ll\frac{1}{n^{2}}\sum_{u=1}^{k_{i}}\sum_{j\in I_{u}^{\prime}}^{{}}%
\Vert{\mathbb{E}}(Y_{ij}|{\mathcal{F}}_{i,u-1}^{Y})\Vert_{2}^{2}\ll\frac
{1}{n^{2}}\sum_{u=1}^{k_{i}}\sum_{j\in I_{u}^{\prime}}^{{}}\Vert{\mathbb{E}%
}(X_{ij}|{\mathcal{F}}_{i,u-1}^{X})\Vert_{2}^{2}\,.
\]
Since ${\mathcal{F}}_{i,u-1}^{X}\subset{\mathcal{G}}_{i,\ell-q}$ for any
$\ell\in\{(u-1)(p+q)+1,\dots,(u-1)(p+q)+p\}$, it follows that
\begin{equation}
T_{3}\ll\frac{1}{n^{2}}\sum_{j=1}^{i}\Vert{\mathbb{E}}(X_{ij}|{\mathcal{G}%
}_{i,j-q})\Vert_{2}^{2}\,\ll\frac{1}{n}\eta_{q}^{2}. \label{T3}%
\end{equation}
To treat $T_{2}$ we proceed as in the proof of relation \eqref{T3}, and infer
that
\begin{equation}
T_{2}\ll\frac{1}{n^{2}}\sum_{j=1}^{i}\Vert{\mathbb{E}}({\widetilde{X}}%
_{ij}|{\mathcal{G}}_{i,j-q})\Vert_{2}^{2}\ll\frac{1}{n}\eta_{q}^{2}+\frac
{1}{n^{2}}\sum_{j=1}^{i}\Vert X_{ij}^{2}I(|X_{ij}|>A)\Vert_{1}\,. \label{T2}%
\end{equation}
We handle now the term $T_{1}$ in \eqref{TErowsinterterme1}. Using the
notation \eqref{defCiu} and the fact that $\mathbf{Y}_{n}$ has the same
covariance structure as $\mathbf{X}_{n}$, we start by rewriting $T_{1}$ as
follows:
\begin{align}
T_{1}  &  =\sum_{u=1}^{k_{i}}\sum_{j,\ell\in I_{u}^{\prime}}^{{}%
}\big |{\mathbb{E}}\big (\big ({\mathbb{E}}({\widetilde{X}}_{ij}{\widetilde
{X}}_{i\ell}|{\mathcal{F}}_{i,u-1}^{X})-{\mathbb{E}}({X}_{ij}{X}_{i\ell
})\big )\partial_{ij}\partial_{i\ell}s(C_{i,u}%
)\big )\big |\nonumber\label{defT1}\\
&  =\sum_{u=1}^{k_{i}}\sum_{j,\ell\in I_{u}^{\prime}}\big |{\mathbb{E}%
}\big (\big ({\widetilde{X}}_{ij}{\widetilde{X}}_{i\ell}-{\mathbb{E}}({X}%
_{ij}{X}_{i\ell})\big )\partial_{ij}\partial_{i\ell}s(C_{i,u})\big )\big |\,,
\end{align}
where for the second equality we used the fact that $\partial_{ij}%
\partial_{i\ell}s(C_{i,u})$ is measurable with respect to ${\mathcal{H}}%
_{i,u}$ defined by \eqref{definitionHX} and that $\sigma
\big ((X_{i,(u-1)(p+q)+j})_{1\leq j\leq p}\big )\vee{\mathcal{F}}_{i,u-1}^{X}$
is independent of
\[
\sigma\big ((L_{j}^{X})_{1\leq j\leq i-1},(L_{k}^{Y})_{i+1\leq k\leq
n}\big )\vee\sigma\big ({D}_{i,u+1}^{\ast},\dots,{D}_{ik_{i}}^{\ast}\big ).
\]

To treat the summands in \eqref{defT1}, we further weaken the dependence by
suppressing some variables in $C_{i,u}$ which are \textquotedblright
close\textquotedblright\ to ${\widetilde{X}}_{ij}{\widetilde{X}}_{i\ell}$. Let
$a$ be a positive integer fixed for the moment. Then, setting,
\[
C_{i,u}^{(a)}=\big ({\widetilde{D}}_{[1,i-1]}^{X},{\widetilde{D}}_{i1}%
^{\prime},\dots,{\widetilde{D}}_{i,u-a}^{\prime},{\mathbf{0}},{D^{\prime}%
}_{i,u+1}^{\ast},\dots,{D^{\prime}}_{i,u_{k_{i}}}^{\ast},{D}_{[i+1,n]}%
^{Y}\big )\,\text{ if $u\geq a+1$}\,,
\]
and
\[
C_{i,u}^{(a)}=\big ({\widetilde{D}}_{[1,i-1]}^{X},{\mathbf{0}},{D^{\prime}%
}_{i,u+1}^{\ast},\dots,{D^{\prime}}_{i,u_{k_{i}}}^{\ast},{D}_{[i+1,n]}%
^{Y}\big )\,\text{ if $1\leq u\leq a$}\,,
\]
we write
\begin{equation}
\big |{\mathbb{E}}\big (\big ({\widetilde{X}}_{ij}{\widetilde{X}}_{i\ell
}-{\mathbb{E}}({X}_{ij}{X}_{i\ell})\big )\partial_{ij}\partial_{i\ell
}s(C_{i,u})\big )\big |\leq I_{1}+I_{2}\,. \label{decT1}%
\end{equation}
where
\[
I_{1}=\big |{\mathbb{E}}\big (\big ({\widetilde{X}}_{ij}{\widetilde{X}}%
_{i\ell}-{\mathbb{E}}({X}_{ij}{X}_{i\ell})\big )\partial_{ij}\partial_{i\ell
}\big (s_{n}(C_{i,u})-s_{n}(C_{i,u}^{(a)})\big)\big )\big |
\]
and
\[
I_{2}=\big |{\mathbb{E}}\big (\big ({\widetilde{X}}_{ij}{\widetilde{X}}%
_{i\ell}-{\mathbb{E}}({X}_{ij}{X}_{i\ell})\big )\partial_{ij}\partial_{i\ell
}s(C_{i,u}^{(a)})\big )\big |\,.
\]

By using the multivariate Taylor expansion of first order for $\partial
_{ij}\partial_{i\ell}s$, taking into account the definitions of $C_{i,u}$ and
$C_{i,u}^{(a)}$ and then by using \eqref{boundsd}, we derive, after simple
computations, that
\begin{equation}
I_{1}\ll\frac{1}{n^{5/2}}\sum_{v=2}^{a+1}\sum_{r\in I_{v}^{\prime}}%
\Vert\big ({\widetilde{X}}_{ij}{\widetilde{X}}_{i\ell}-{\mathbb{E}}({X}%
_{ij}{X}_{i\ell})\big )\big({\widetilde{X}}_{ir}-{\mathbb{E}}({\widetilde{X}%
}_{ir}|{\mathcal{F}}_{u-v}^{X}\big )\Vert_{1}\ll\frac{1}{n^{5/2}}%
(Aap)\sigma^{2}\,. \label{T11}%
\end{equation}
Next, using \eqref{boundsd} again and the definition of the conditional
expectation, we infer that
\[
I_{2}\ll\frac{1}{n^{2}}\Vert{\mathbb{E}}\big ({\widetilde{X}}_{ij}%
{\widetilde{X}}_{i\ell}|\sigma(C_{i,u}^{(a)})\big )-{\mathbb{E}}({X}_{ij}%
{X}_{i\ell})\Vert_{1}\,.
\]
Notice now that, since $\mathbf{X}_{n}$ and $\mathbf{Y}_{n}$ are assumed to be
independent and since the rows of ${\mathbf{X}}_{n}$ are independent,
${\mathbb{E}}\big ({\widetilde{X}}_{ij}{\widetilde{X}}_{i\ell}|\sigma
(C_{i,u}^{(a)})\big )={\mathbb{E}}\big ({\widetilde{X}}_{ij}{\widetilde{X}%
}_{i\ell}|{\mathcal{F}}_{i,u-a}^{X}\big )$. Therefore, after simple
computations based on the definition of ${\widetilde{X}}_{ij}$ and on the fact
that $A\Vert X_{ij}I(|X_{ij}|>A)\Vert_{1}\leq$ $\Vert X_{ij}^{2}%
I(|X_{ij}|>A)\Vert_{1}$, we obtain
\begin{equation}
I_{2}\ll\frac{1}{n^{2}}\Vert{\mathbb{E}}\big ({X}_{ij}{X}_{i\ell}%
|{\mathcal{F}}_{i,u-a}^{X}\big )-{\mathbb{E}}({X}_{ij}{X}_{i\ell})\Vert
_{1}+\frac{1}{n^{2}}\Vert X_{ij}I(|X_{ij}|>A)\Vert_{2}\Vert X_{i\ell
}I(|X_{i\ell}|>A)\Vert_{2}\,. \label{T12}%
\end{equation}
Starting from \eqref{defT1} and taking into account \eqref{decT1}, \eqref{T11}
and \eqref{T12}, we get
\begin{equation}
T_{1}\ll\frac{1}{n^{3/2}}(Aap^{2})\sigma^{2}+\frac{p}{n^{2}}\sum_{j=1}%
^{i}\Vert X_{ij}^{2}I(|X_{ij}|>A)\Vert_{1}+\frac{1}{n^{2}}k_{i}p^{2}%
\gamma_{aq}\,. \label{T1}%
\end{equation}
So, overall, starting now from the inequality \eqref{TErowsinter}, taking into
account \eqref{TErowsinterterme1}, \eqref{T3}, \eqref{T2} and \eqref{T1}, and
summing over $i$, we obtain that
\begin{equation}
\big |\mathbb{E}s_{n}({\widetilde{D}}^{X})-\mathbb{E}s_{n}(D^{Y}%
)\big |\ll\frac{1}{n^{1/2}}p^{2}\sigma^{2}(A+aA+\sigma)+pL(A)+\eta_{q}%
^{2}+p\gamma_{aq}\,. \label{step3}%
\end{equation}
\medskip

\noindent\textbf{Step 4: End of the proof.}

Starting from \eqref{mart decompose}, taking $A=\varepsilon\sqrt{n}$ and
considering the upper bound \eqref{step3}, we get
\begin{multline*}
\big |\mathbb{E}S^{\mathbb{X}_{n}}(z)-\mathbb{E}S^{\mathbb{Y}_{n}}(z)\big |\ll
p^{2}\sigma^{2}(\varepsilon+a\varepsilon+\frac{1}{n^{1/2}}\sigma
)+pL(\varepsilon\sqrt{n})+\eta_{q}^{2}+p\gamma_{aq}\\
+\Big (\frac{q}{p}+\frac{q+p}{n}\Big )^{1/2}\sigma+L^{1/2}(\varepsilon\sqrt
{n})+\eta_{q}\,.
\end{multline*}
Therefore, when $n\rightarrow\infty,$ we obtain for all $p,q,a,$ and
$\varepsilon$,
\[
\limsup_{n\rightarrow\infty}\big |\mathbb{E}S^{\mathbb{X}_{n}}(z)-\mathbb{E}%
S^{\mathbb{Y}_{n}}(z)\big |\ll p^{2}\sigma^{2}(\varepsilon+a\varepsilon
)+\eta_{q}^{2}+\eta_{q}+p\gamma_{aq}+(q/p)^{1/2}\,\sigma.
\]
Now we let $\varepsilon\rightarrow0$ and obtain%
\[
\limsup_{n\rightarrow\infty}\big |\mathbb{E}S^{\mathbb{X}_{n}}(z)-\mathbb{E}%
S^{\mathbb{Y}_{n}}(z)\big |\ll\eta_{q}^{2}+\eta_{q}+p\gamma_{aq}%
+(q/p)^{1/2}\,\sigma.
\]
Then we let $a\rightarrow\infty,$ and, by our hypotheses, for any $p$ and $q$
we obtain
\[
\limsup_{n\rightarrow\infty}\big |\mathbb{E}S^{\mathbb{X}_{n}}(z)-\mathbb{E}%
S^{\mathbb{Y}_{n}}(z)\big |\ll\eta_{q}^{2}+\eta_{q}+(q/p)^{1/2}\,\sigma.
\]
Now we can let $p$ and $q$ tend to $\infty$ in such a way $q/p\rightarrow0$ to
obtain the desired result. \ $\ \ \ \lozenge$

\subsection{Proof of Corollary \ref{stationary}}

By the reverse martingale convergence theorem and condition \eqref{regular 1},
we get that $\lim_{n\rightarrow\infty}\mathbb{E}(X_{0}|\mathcal{G}_{-n}%
)$\newline$=\mathbb{E}(X_{0}|{\mathcal{G}}_{-\infty})=0$ a.s. So, since
$X_{0}$ belongs to ${\mathbb{L}}^{2}$, this last convergence implies that
condition (\ref{mart-K-indrow}) holds. We prove now that under the conditions
of the corollary, condition \eqref{mixing-indrow} is satisfied. Note first
that, by stationarity, this latter condition reads as
\begin{equation}
\sup_{u}\Vert\mathbb{E}(X_{0}X_{u}|\mathcal{G}_{-n})-\mathbb{E}(X_{0}%
X_{u})\Vert_{1}\rightarrow0\text{ as }n\rightarrow\infty
\,.\label{cond2mixindrows}%
\end{equation}
To prove that \eqref{cond2mixindrows} holds we shall prove that
\begin{equation}
\lim_{p\rightarrow\infty}\limsup_{n\rightarrow\infty}\sup_{u\geq p+1}%
\Vert\mathbb{E}(X_{0}X_{u}|\mathcal{G}_{-n})-\mathbb{E}(X_{0}X_{u})\Vert
_{1}=0\,,\label{cond2mixindrows1}%
\end{equation}
and that
\begin{equation}
\lim_{p\rightarrow\infty}\limsup_{n\rightarrow\infty}\max_{1\leq u\leq p}%
\Vert\mathbb{E}(X_{0}X_{u}|\mathcal{G}_{-n})-\mathbb{E}(X_{0}X_{u})\Vert
_{1}=0\,.\label{cond2mixindrows2}%
\end{equation}
To prove \eqref{cond2mixindrows1}, we note that
\begin{align*}
\sup_{u\geq p+1}\Vert\mathbb{E}(X_{0}X_{u}|\mathcal{G}_{-n})-\mathbb{E}%
(X_{0}X_{u})\Vert_{1} &  \leq\sup_{u\geq p+1}\Vert\mathbb{E}(X_{0}%
X_{u}|\mathcal{G}_{0})-\mathbb{E}(X_{0}X_{u})\Vert_{1}\\
&  =\sup_{u\geq p+1}\Vert X_{0}\mathbb{E(}X_{u}|\mathcal{G}_{0})-\mathbb{E}%
(X_{0}X_{u})\Vert_{1}\\
&  \leq2\Vert X_{0}\Vert_{2}\cdot\sup_{u\geq p+1}\Vert\mathbb{E(}%
X_{u}|\mathcal{G}_{0})\Vert_{2}\leq2\Vert X_{0}\Vert_{2}\cdot\Vert
\mathbb{E(}X_{0}|\mathcal{G}_{-p})\Vert_{2}\,.
\end{align*}
This shows that \eqref{cond2mixindrows1} holds since (\ref{mart-K-indrow})
does under \eqref{regular 1}. We turn now to the proof of
\eqref{cond2mixindrows2}. By the reverse martingale convergence theorem
\begin{gather*}
\lim_{n\rightarrow\infty}\max_{1\leq u\leq p}\Vert\mathbb{E}(X_{0}%
X_{u}|\mathcal{G}_{-n})-\mathbb{E}(X_{0}X_{u})\Vert_{1}=\max_{1\leq u\leq
p}\lim_{n\rightarrow\infty}\Vert\mathbb{E}(X_{0}X_{u}|\mathcal{G}%
_{-n})-\mathbb{E}(X_{0}X_{u})\Vert_{1}\\
=\sup_{1\leq u\leq p}\Vert\mathbb{E}(X_{0}X_{u}|\mathcal{G}_{-{\infty}%
})-\mathbb{E}(X_{0}X_{u})\Vert_{1}\,,
\end{gather*}
which is equal to zero by condition (\ref{regular 2}). This ends the proof of
\eqref{cond2mixindrows2} and then of the corollary. $\lozenge$

\subsection{Proof of Theorem \ref{Gram}}

It is well-known that for deriving the limiting spectral distribution of
$\mathbb{B}_{N}$ it is enough to study the Stieltjes transform of the
following symmetric matrix of order $n=N+p$:
\[
\mathbb{X}_{n}=\frac{1}{\sqrt{N}}\left(
\begin{array}
[c]{cc}%
\mathbf{0}_{p,p} & {\mathcal{X}}_{N,p}^{T}\\
{\mathcal{X}}_{N,p} & \mathbf{0}_{N,N}%
\end{array}
\right)  \,.
\]
Indeed the eigenvalues of $\mathbb{X}_{n}^{2}$ are the eigenvalues of
$N^{-1}{\mathcal{X}}_{N,p}^{T}{\mathcal{X}}_{N,p}$ together with the
eigenvalues of $N^{-1}{\mathcal{X}}_{N,p}{\mathcal{X}}_{N,p}^{T}$. Since these
two latter matrices have the same nonzero eigenvalues, the following relation
holds: for any $z\in{\mathbb{C}}^{+}$, $S_{\mathbb{B}_{N}}(z)=z^{-1/2}\frac
{n}{2p}S_{\mathbb{X}_{n}}(z^{1/2})+\frac{N-p}{2pz}$ (see, for instance, page
549 in Rashidi Far \textit{et al.} \cite{ROBS} for additional arguments
leading to the relation above. Obviously a similar equation holds for the Gram
random matrix $\mathbb{H}_{N}$ associated with $(Y_{\mathbf{u}})_{\mathbf{u}%
\in{\mathbb{Z}}^{2}}$, namely: $S_{\mathbb{H}_{N}}(z)=z^{-1/2}\frac{n}%
{2p}S_{\mathbb{Y}_{n}}(z^{1/2})+\frac{N-p}{2pz}$, where $\mathbb{Y}_{n}$ is
defined as $\mathbb{X}_{n}$ but with $X_{\mathbf{u}}$ replaced by
$Y_{\mathbf{u}}$. Therefore, in order to prove the theorem, it suffices to
show that, for any $z\in{\mathbb{C}}^{+}$,
\begin{equation}
\lim_{N\rightarrow\infty}\big |S_{\mathbb{X}_{n}}(z)-{\mathbb{E}%
}(S_{\mathbb{Y}_{n}}(z))\big |=0\,\text{ a.s.} \label{butgram}%
\end{equation}
Note now that $\mathbb{X}_{n}:=n^{-1/2}[x_{ij}^{(n)}]_{i,j=1}^{n}$ where
$x_{ij}^{(n)}=\sqrt{\frac{n}{N}}X_{i-p,j}{\mathbf{1}}_{i\geq p+1}{\mathbf{1}%
}_{1\leq j\leq p}$ if $1\leq j\leq i\leq n$, and $x_{ij}^{(n)}=x_{ji}^{(n)}$
if $1\leq i<j\leq n$. Similarly we can write $\mathbb{Y}_{n}:=n^{-1/2}%
[y_{ij}^{(n)}]_{i,j=1}^{n}$ where the $y_{ij}^{(n)}$'s are defined as the
$x_{ij}^{(n)}$'s but with $X_{i-p,j}$ replaced by $Y_{i-p,j}$. The theorem
then follows by applying Remark \ref{remtriangular} of Theorem
\ref{mainindependentrows} to the matrices $\mathbb{X}_{n}$ and $\mathbb{Y}%
_{n}$ defined above. \ $\ \ \ \lozenge$

\subsection{Proof of Theorem \ref{corgram}}

According to Theorem \ref{Gram} and Theorem B.9. in Bai and Silverstein
(2010), the proof of Theorem \ref{corgram} is reduced to show that, for any
$z\in{\mathbb{C}}^{+}$
\begin{equation}
\lim_{N\rightarrow\infty}{\mathbb{E}}(S^{\mathbb{H}_{N}}(z))=S(z),
\label{conv1corgram}%
\end{equation}
where $\mathbb{H}_{N}$ is the Gram matrix associated with a Gaussian random
field $(Y_{\mathbf{u}})_{\mathbf{u}\in{\mathbb{Z}}}$ having the same
covariance structure as $(X_{\mathbf{u}})_{\mathbf{u}\in{\mathbb{Z}}}$ and
$S(z)$ is a Stieltjes transform of a measure $F$ and satisfies equation
\eqref{equationlimit}. To prove the convergence above, we shall proceed in two
steps. In the first step we shall prove that \eqref{conv1corgram} holds under
the additional assumption that the spectral density of $(X_{k})_{k\in
{\mathbb{Z}}}$ is square integrable. The proof of this particular case is
facilitated by the fact that a square integrable spectral density allows us to
use the celebrated Szeg\"{o}-Trotter theorem for Toeplitz matrices. This
assumption will be removed in a second step, where we approximate the spectral
density by a square integrable one and then extend the characterization of the limit.

\medskip

\noindent\textbf{Step 1. Proof of \eqref{conv1corgram} when the spectral
density is square integrable.}

We shall apply Theorem 1.1 in Silverstein (1995). Consider $N$ independent
copies $(g_{ij})_{j\in{\mathbb{Z}}}$, $i=1,\dots,N$ of a sequence
$(g_{k})_{k\in{\mathbb{Z}}}$ of i.i.d. standard normal random variables. Set
\[
\Gamma_{p}:=\left(
\begin{array}
[c]{cccc}%
c_{0} & c_{1} & \cdots & c_{p-1}\\
c_{1} & c_{0} &  & c_{p-2}\\
\vdots & \vdots & \vdots & \vdots\\
c_{p-1} & c_{p-2} & \cdots & c_{0}%
\end{array}
\right)  \ \text{ where}\ c_{k}=\mathrm{Cov}(X_{0},X_{k})\,.
\]
Using the stationarity of the Gaussian process $(Y_{\mathbf{u}})_{{\mathbf{u}%
}\in{\mathbb{Z}}^{2}}$, we can easily verify that the random vector
$((Y_{1j})_{1\leq j\leq p},\dots,(Y_{Nj})_{1\leq j\leq p})$ has the same
distribution as $(\mathbf{g}_{1}\Gamma_{p}^{1/2},\dots,\mathbf{g}_{N}%
\Gamma_{p}^{1/2})$ where for any $i\in\{1,\dots,N\}$, $\mathbf{g}_{i}%
=(g_{ij})_{1\leq j\leq p}$ and $\Gamma_{p}^{1/2}$ is the symmetric
non-negative square root of $\Gamma_{p}$. Therefore, for any $z\in{\mathbb{C}%
}^{+}$,
\[
{\mathbb{E}}(S^{\mathbb{H}_{N}}(z))={\mathbb{E}}(S^{\Gamma_{p}^{1/2}%
\mathbb{G}_{N}\Gamma_{p}^{1/2}}(z))\,,
\]
where $\mathbb{G}_{N}=\frac{1}{N}{\mathcal{G}}_{N,p}^{T}{\mathcal{G}}_{N,p}$
with ${\mathcal{G}}_{N,p}=(g_{ij})_{1\leq i\leq N,1\leq j\leq p}$. Hence,
according to Theorem 1.1 in Silverstein (1995), if $p/N\rightarrow
c\in(0,\infty)$ and
\begin{equation}
F^{\Gamma_{p}}\text{ converges to a probability distribution $H$ as
$p\rightarrow\infty$}, \label{convdistgammaN}%
\end{equation}
then there is a nonrandom probability distribution $F$ such that
\begin{equation}
\label{resultatsilverstein}d(F^{{\mathbb{H}}_{N}},F)\rightarrow0 \, \text{
a.s.}%
\end{equation}
Furthermore, the Stieltjes transform $S=S(z),$ $z\in\mathbb{C}^{+}$, of $F$
satisfies the equation
\[
S=\int\frac{1}{x(1-c-czS)-z}dH(x).\
\]
Setting ${\underline{S}}:=-(1-c)/z+cS$, this last equation becomes
\begin{equation}
z=-\frac{1}{{\underline{S}}}+c\int\frac{x}{1+x{\underline{S}}}dH(x).
\label{equation}%
\end{equation}
We mention that ${\underline{S}}$ is also a Stieltjes transform (see relation
(1.3) in {\cite{Silver} or \cite{GH}), so }$\operatorname{Im}{\underline{S}%
>0}$ for $z\in\mathbb{C}^{+}.$

Note now that, since the spectral density $f$ is assumed to be square
integrable, by Parseval's identity we have that $\sum_{k\in{\mathbb{Z}}}%
c_{k}^{2}<\infty$. Therefore by a version of the Szeg\"{o}'s theorem for
Toeplitz forms (see page 72 of Trotter (1984)), the convergence
\eqref{convdistgammaN} holds and we have, for any $\varphi$ which is
continuous and bounded,
\[
\int\varphi(x)dH(x)=\frac{1}{2\pi}\int_{-\pi}^{\pi}\varphi(2\pi f(\lambda
))d\lambda\,.
\]
Since the function $\varphi(x):=x/(1+x{\underline{S}})$ is continuous and
bounded by $1$/$\operatorname{Im}{\underline{S}}$, the relation
(\ref{equation}) can be rewritten as (\ref{equationlimit}). To end the proof
of \eqref{conv1corgram} when the spectral density is assumed to be square
integrable, it suffices to notice that \eqref{resultatsilverstein} implies
that $\lim_{N\rightarrow\infty}S^{\mathbb{H}_{N}}(z)=S(z)$ a.s. which in turn
entails \eqref{conv1corgram} since the Stieltjes transforms are bounded.

\medskip

\noindent\textbf{Step 2. Proof of \eqref{conv1corgram} when the spectral
density is not necessarily square integrable.}

To remove the assumption on the square integrability of the spectral density,
we shall truncate the spectral density, then define a Gaussian process with
the help of the truncated spectral density. Next, we use the limit of the
empirical eigenvalue distribution for this truncated process to approximate
and then characterize the limit of $F^{\mathbb{H}_{n}}$.

In the rest of the proof, $(\xi_{{\mathbf{u}}})_{{\mathbf{u}}\in{\mathbb{Z}%
}^{2}}$ is a sequence of i.i.d. standard normal real-valued random variables.
According to Proposition \ref{prop GaussLinear}, there is no loss of
generality by assuming from now on that $(Y_{\mathbf{u}})_{{\mathbf{u}}%
\in{\mathbb{Z}}}$ has the following linear representation: for any $k,\ell$ in
${\mathbb{Z}}$,
\begin{equation}
Y_{k\ell}=\sum_{j\in{\mathbb{Z}}}a_{j}\xi_{k,\ell-j}\,\text{ with }%
\,a_{k}=\frac{1}{\sqrt{2\pi}}\int_{-\pi}^{\pi}e^{{\mathrm{i}}kx}\sqrt
{f(x)}dx\,.\label{defakYkf}%
\end{equation}
For a fixed positive real $b$, we define another centered real-valued Gaussian
random field $(Z_{\mathbf{u}}^{b})_{{\mathbf{u}}\in{\mathbb{Z}}^{2}}$ with the
help of the function
\[
f_{b}=f\wedge b\,.
\]
Note that since $f$ is a nonnegative, even and integrable function on
$[-\pi,\pi]$, so is $f_{b}$. Then $f_{b}$ is also the spectral density on
$[-\pi,\pi]$ of a ${\mathbb{L}}^{2}$-stationary process. Therefore, according
to Proposition \ref{prop GaussLinear}, if we set, for any $k,\ell$ in
${\mathbb{Z}}$,
\begin{equation}
{\tilde{a}}_{k}=\frac{1}{\sqrt{2\pi}}\int_{-\pi}^{\pi}e^{{\mathrm{i}}kx}%
\sqrt{f_{b}(x)}dx\,\text{ and }\,Z_{k\ell}^{b}=\sum_{j\in{\mathbb{Z}}}%
{\tilde{a}}_{j}\xi_{k,\ell-j}\,,\label{defakZkf}%
\end{equation}
$(Z_{\mathbf{u}}^{b})_{{\mathbf{u}}\in{\mathbb{Z}}^{2}}$ is a centered
real-valued stationary Gaussian random field. In addition, for any fixed
integer $k$, $(Z_{k\ell}^{b})_{\ell\in{\mathbb{Z}}}$ admits $f_{b}$ as
spectral density on $[-\pi,\pi]$. Let $\mathbb{H}_{N}^{b}$ be the Gram matrix
associated with $(Z_{\mathbf{u}}^{b})_{{\mathbf{u}}\in{\mathbb{Z}}^{2}}$.
Since $f_{b}$ is bounded, it is in particular square integrable. Then, by the
Step 1 of the proof, we conclude that there is a nonrandom distribution
function $F^{b}$ such that
\begin{equation}
\lim_{N\rightarrow\infty}d(F^{\mathbb{H}_{N}^{b}}\ ,F^{b}\ )=0\text{
a.s.}\label{conv3corgram}%
\end{equation}
On another hand, by using Lemma \ref{lmagotze2} together with Cauchy-Schwarz's
inequality, we infer that
\[
{\mathbb{E}}d^{2}(F^{\mathbb{H}_{N}}\ ,F^{b}\ )\ \ll\frac{1}{Np}%
\Big \Vert\sum_{i=1}^{N}\sum_{j=1}^{p}(Y_{ij}^{2}+(Z_{ij}^{b})^{2}%
)\Big \Vert_{1}^{{1/2}}\Big \Vert\sum_{i=1}^{N}\sum_{j=1}^{p}(Y_{ij}%
-Z_{ij}^{b})^{2}\Big \Vert_{1}^{{1/2}}\,.
\]
$\ $Since ${\mathbb{E}}(Y_{ij}^{2})=\sum_{k\in{\mathbb{Z}}}a_{k}^{2}$,
${\mathbb{E}}((Z_{ij}^{b})^{2})=\sum_{k\in{\mathbb{Z}}}{\tilde{a}}_{k}^{2}$
and ${\mathbb{E}}((Y_{ij}^{{}}-Z_{ij}^{b})^{2})=\sum_{k\in{\mathbb{Z}}}%
(a_{k}-{\tilde{a}}_{k})^{2}$, by using Parseval's identity, it follows that
\begin{multline*}
{\mathbb{E}}d^{2}(F^{\mathbb{H}_{N}},F^{b})\ll\Big (\int_{-\pi}^{\pi
}f(x)dx+\int_{-\pi}^{\pi}f_{b}(x)dx\Big )^{1/2}\Big (\int_{-\pi}^{\pi}%
(f^{1/2}(x)-f_{b}^{1/2}(x))^{2}dx\Big )^{1/2}\\
\ll\Big (\int_{-\pi}^{\pi}f(x)dx\Big )^{1/2}\Big (\int_{-\pi}^{\pi
}fI(f>b)(x)dx\Big )^{1/2}\,.
\end{multline*}
Therefore, by the Lebesgue dominated convergence theorem%
\begin{equation}
\lim_{b\rightarrow\infty}\limsup_{N\rightarrow\infty}{\mathbb{E}}%
d^{2}(F^{\mathbb{H}_{N}},F^{b})=0\,.\label{conv4corgram}%
\end{equation}
Since for any positive reals $b$ and $b^{\prime}$ we have
\[
d^{2}(F^{b^{\prime}},F^{b})\leq2{\mathbb{E}}d^{2}(F^{\mathbb{H}_{N}%
},F^{b^{\prime}})+2{\mathbb{E}}d^{2}(F^{\mathbb{H}_{N}},F^{b})\,,
\]
(\ref{conv4corgram}) implies that $F^{b}$ is Cauchy. Taking into account that
the space of distribution functions endowed with L\'{e}vy metric is complete,
we conclude that there is a nonrandom distribution function $F$ such that
$\lim_{b\rightarrow\infty}d(F^{b},F)=0$ which, combined with
(\ref{conv4corgram}), also gives $\lim_{N\rightarrow\infty}{\mathbb{E}%
}d(F^{\mathbb{H}_{N}},F\ )=0.$ If we denote by $S$ the Stieltjes transform of
$F$ and by $S^{b}$ the Stieltjes transform of $F^{b}$, by the continuity
theorem (see for instance Theorem B.9 in \cite{BS}), we obtain, for any
$z\in\mathbb{C}^{+}$, the convergence of $S^{b}(z)$ to $S(z)$ and the
convergence in probability of $S^{\mathbb{H}_{N}}(z)$ to $S(z)$. Since the
Stieltjes transforms are bounded, we also have $\lim_{N\rightarrow\infty
}{\mathbb{E}}(S^{\mathbb{H}_{N}}(z))=S(z),$ which completes the proof of the
convergence (\ref{conv1corgram}).

We shall prove now that $S(z)$ satisfies (\ref{equationlimit}). We start from
the equation satisfied by $S^{b}$ which was found in Step 1, namely
\begin{equation}
z=-\frac{1}{{\underline{S}}^{b}}+c\int_{-\pi}^{\pi}\frac{f_{b}(x)}{1+2\pi
f_{b}(x){\underline{S}}^{b}}dx\,,\label{equation2}%
\end{equation}
with ${\underline{S}}^{b}:=-(1-c)/z+cS^{b}$. We note at this point that, and
for any $z$ in $\mathbb{C}^{+}$, we also have ${\underline{S}}(z)=\lim
_{b\rightarrow\infty}{\underline{S}}^{b}(z).$ $\ $It follows that
${\underline{S}=}-(1-c)/z+cS,$ where ${\underline{S}}$ is also a Stieltjes
transform, implying that $\operatorname{Im}({\underline{S}})(z)>0$. Therefore,
for almost all $x$ in $[-\pi,\pi]$,
\[
\lim_{b\rightarrow\infty}\frac{f_{b}(x)}{1+2\pi f_{b}(x){\underline{S}}^{b}%
}=\frac{f(x)}{1+2\pi f(x){\underline{S}}}\,.
\]
Also, for all $b$ sufficiently large,
\[
\Big |\frac{f_{b}(x)}{1+2\pi f_{b}(x){\underline{S}}^{b}}\Big |\leq\frac
{1}{2\pi\operatorname{Im}({\underline{S}}^{b})}\leq\frac{1}{\operatorname{Im}%
({\underline{S}})}\,.
\]
By the Lebesgue dominated convergence theorem, by passing to the limit when
$b\rightarrow\infty$ in (\ref{equation2}) we obtain that ${\underline{S}}$
satisfies equation (\ref{equationlimit}).$\ \ \ \ \ \lozenge$

\end{document}